\documentclass[11pt]{article}
\usepackage{ytableau}
\usepackage{bm}
\usepackage{array}
\usepackage[margin=1in]{geometry}
\usepackage{mathrsfs}
\usepackage{stmaryrd}
\usepackage{amsmath,amssymb,amscd}
\usepackage{shadow}
\usepackage{epsfig,epstopdf}
\usepackage{graphicx}
\usepackage{color}
\usepackage{multirow}
\usepackage{booktabs} %调整表格线与上下内容的间隔
\usepackage{enumerate}
\usepackage{longtable}
\usepackage{pstricks,multido}
\usepackage[colorlinks,urlcolor=purple,linkcolor=red,anchorcolor=green,citecolor=blue]{hyperref}
\usepackage{authblk}
\allowdisplaybreaks
\newtheorem{thm}{Theorem}[section]

\newtheorem{cor}[thm]{Corollary}

\newtheorem{example}{Example}[section]

\numberwithin{figure}{section}
\numberwithin{table}{section}
\def\qed{\hfill \rule{4pt}{7pt}}
\def\pf{\noindent {\it{Proof.} \hskip 2pt}}

\def\gen{{\rm Gen}}
\def\erf{{\rm erf}}
\numberwithin{equation}{section}
\begin{document}

\begin{center}
{\Large\bf Joint Distributions  of
Permutation Statistics and \\[6pt]  the
 Parabolic Cylinder Functions}
\end{center}

\begin{center}

Amy M. Fu$^{1}$
and Frank Z.K. Li$^{2}$

$^{1}$School of Mathematics,\\ Shanghai University of Finance and Economics,\\
 Shanghai 200433, P.R. China\\[6pt]

$^{2}$Center for Combinatorics, LPMC, \\
Nankai University, Tianjin 300071, P.R. China\\[6pt]

Emails:
{\tt $^{1}$fu.mei@mail.shufe.edu.cn},
{\tt $^{2}$zhkli@mail.nankai.edu.cn}

\end{center}

%\maketitle

\noindent{\bf Abstract}:
%\begin{abstract}
%The concept of the grammatical labeling was introduced
%by Chen and Fu to build the connections between context-free grammars and the combinatorial structures.
In this paper, we introduce a
context-free grammar $G\colon x \rightarrow  xy,\, y \rightarrow  zu,\,
       z \rightarrow  zw,\,  w \rightarrow  xv,\,
       u \rightarrow xyz^{-1}v,\,
       v \rightarrow x^{-1}zwu$
over the variable set $V=\{x,y,z,w,u,v\}$.
We use this grammar to study  joint distributions  of
several permutation statistics related to descents, rises,
peaks and valleys. By considering the pattern  of an exterior peak, we introduce the exterior peaks of pattern 132 and
of pattern 231. Similarly, peaks can also be classified
according to their patterns.
Let $D$ be the formal derivative operator with respect to the grammar $G$.
By using a grammatical labeling,
we show  that $D^n(z)$ is the generating function
of the number of permutations on $[n]=\{1,2,\ldots,n\}$
with given numbers of exterior peaks of pattern 132
and of pattern 231, and proper double descents.
By solving a cylinder differential  equation, we obtain
an  explicit formula of the generating function of $D^n(z)$,
which can be viewed as a unification of the results
of Elizalde-Noy, Barry, Basset, Fu and Gessel.
Specializations lead to
the joint distributions of certain consecutive patterns in permutations,
 as studied by Elizalde-Noy and Kitaev.
By a different labeling with respect to the same grammar $G$,
we  derive  the joint distribution of peaks of pattern 132 and of pattern 231,
double descents and  double rises, with the
 generating function also  expressed by the parabolic cylinder functions. This formula serves as a refinement of the work of Carlitz-Scoville.
Furthermore, we obtain  the joint distribution of  exterior peaks of pattern 132
and  of pattern 231 over alternating permutations.
%\end{abstract}



\section{Introduction}

In this paper, we introduce a refinement of the
exterior peaks based on their patterns.
More precisely,
an exterior peak of a permutation is either of pattern 132 or
231.
Similarly, we may define the types of a peak of a permutation.
By giving a context-free grammar, we study the joint distribution of  exterior peaks of pattern 132 and of pattern 231,
and proper double descents over permutations on $[n]=\{1, 2, \ldots, n\}$.
In the same manner, we study the joint distribution of  peaks of
pattern 132 and of pattern 231, double descents and double rises on permutations.
It turns out that the generating functions for the joint distributions are in connection with
the parabolic cylinder functions.

The notion of exterior peaks was introduced by Aguiar, Bergeron and Nymanin \cite{ABN-2004}.
Let $\pi=\pi_1\pi_2\cdots\pi_n$
be a permutation on $[n]$ and $\pi_0=0$,
then for the index $1\leq i\leq n-1$,
we call $i$ an exterior peak
if $\pi_{i-1}<\pi_i>\pi_{i+1}$.
For example, the permutation $534621$ has two exterior peaks 1 and 4.
For $n\geq0$, denote by $T(n,k)$ the number of permutations on $[n]$ with $k$
exterior peaks and let
$$T_n(x)=\sum_{k\geq0}T(n,k)x^k.$$
Gessel \cite[A008971]{OEIS} obtained  the generating function of $T_n(x)$.

\begin{thm}\emph{\bf(Gessel \cite{OEIS})}\label{Gesselthm}
  We have
  \begin{equation*}
  \sum_{n=0}^\infty T_n(x)\frac{t^n}{n!}=\frac{\sqrt{1-x}}{\sqrt{1-x}\cosh(\sqrt{1-x}t)-\sinh(\sqrt{1-x}t)}.
  \end{equation*}
\end{thm}

Given a permutation $\pi=\pi_1\pi_2\cdots\pi_n$ on $[n]$,
for $2\leq i\leq n-1$,
an index $i$ is a proper double descent if $\pi_{i-1}>\pi_i>\pi_{i+1}$.
As an example, the permutation $653421$ has two proper double descents $2$ and $5$.
For $n\geq0$, let $U(n,k)$ be the number of permutations on $[n]$ with $k$ proper double descents and let
\[U_n(y)=\sum_{k\geq0}U(n,k)y^k.\]
By solving a second-order ordinary differential equation,
Elizalde and Noy \cite{Elizalde-Noy 2003} obtained the generating function of $U_n(y)$.

\begin{thm}\emph{\bf(Elizalde-Noy \cite{Elizalde-Noy 2003})}\label{ENthm}
  We have
  \begin{equation*}
    \sum_{n=0}^\infty U_n(y)\frac{t^n}{n!}=
    \frac{2\sqrt{(y-1)(y+3)}e^{(1-y+\sqrt{(y-1)(y+3)})\frac{t}{2}}}{1+y+\sqrt{(y-1)(y+3)}-(1+y-\sqrt{(y-1)(y+3)})e^{\sqrt{(y-1)(y+3)}t}}.
  \end{equation*}
\end{thm}

Barry \cite{Barry-2014} and Basset \cite{Basset-2014} independently derived the
generating function of permutations on $[n]$ containing no proper double descents.

\begin{thm}\emph{\bf(Elizalde-Noy \cite{Elizalde-Noy 2003}, Barry \cite{Barry-2014}, Basset \cite{Basset-2014})}\label{ENBBthm}
  We have
  \begin{equation*}
    \sum_{n=0}^\infty U(n,0)\frac{t^n}{n!}=
    \frac{\sqrt{3}}{2}\frac{e^{\frac{t}{2}}}{\cos\left(\frac{\sqrt{3}t}{2}+\frac{\pi}{6}\right)}.
  \end{equation*}
\end{thm}

The first named author \cite{Fu-2018} introduced the context-free grammar
\begin{equation}\label{g3}
x\rightarrow xy,\,y\rightarrow xz,\,
z\rightarrow zw,\,w\rightarrow xz,
\end{equation}
and obtained the joint distribution of exterior peaks and
proper double descents over permutations on $[n]$.
For $n\geq0$, let
\begin{equation}\label{fupijk}
  P_n(x,y,z,w)=\sum_{\pi\in\mathfrak{S}_n}x^{{ep}(\pi)}y^{{pdd}(\pi)}z^{{ep}(\pi)+1}w^{n-2ep(\pi)-pdd(\pi)},
\end{equation}
where $\mathfrak{S}_n$ is the set of permutations on $[n]$,
${ep}(\pi)$ and ${pdd}(\pi)$ are the number of exterior peaks and the number of proper double descents in permutation $\pi$, respectively.

\begin{thm}\emph{\bf(Fu \cite{Fu-2018})}\label{futhm}
  We have
  \begin{align*}
    \sum_{n=0}^\infty P_n(x,y,z,w)\frac{t^n}{n!}=
    \frac{2z\sqrt{(y+w)^2-4xz}e^{\frac{t}{2}(w-y+\sqrt{(y+w)^2-4xz})}}{y+w+\sqrt{(y+w)^2-4xz}-(y+w-\sqrt{(y+w)^2-4xz})e^{t\sqrt{(y+w)^2-4xz}}}.
  \end{align*}
\end{thm}

Notice that Theorem \ref{futhm} reduces to Theorem \ref{Gesselthm}, Theorem \ref{ENthm} and Theorem \ref{ENBBthm}.

In this paper,  we give a refinement of the generating function $P_n(x,y,z,w)$
by considering the pattern of an exterior peak.
For a permutation $\pi=\pi_1\pi_2\cdots\pi_n$ and $1\leq i\leq n-1$,
if $i$ is an exterior peak with $\pi_{i-1}<\pi_i>\pi_{i+1}$,
then we say $i$ is an exterior peak of pattern 132
if $\pi_{i-1}<\pi_{i+1}$,
or an exterior peak of pattern 231 if $\pi_{i-1}>\pi_{i+1}$.
For example, the permutation $534621$ has one index 1 as an exterior peak of pattern 132, and one index 4 as an exterior peak of pattern 231.

For $n\geq0$, let $P_n(i,j,k)$ be the number of permutations on $[n]$
with $i$  exterior peaks of pattern 132,
$j$ exterior peaks of pattern 231,
and $k$  proper double descents.
Define
\begin{equation}\label{pnijk}
  P_n( x, y, z, w, u, v)=\sum_{i,j,k=0}^{n}P_n(i,j,k) x^i v^i u^j z^{j+1} y^k w^{n-2i-2j-k},
\end{equation}
or equivalently,
\begin{equation}\label{pnijk2}
  P_n( x, y, z, w, u, v)=\sum_{\pi\in\mathfrak{S}_n}x^{{ep}_1(\pi)}v^{{ep}_1(\pi)} u^{{ep}_2(\pi)} z^{{ep}_2(\pi)+1} y^{{pdd}(\pi)} w^{n-2(ep_1(\pi)+ep_2(\pi))-pdd(\pi)},
\end{equation}
where ${ep}_1(\pi)$ and ${ep}_2(\pi)$ are the numbers of exterior peaks of pattern 132 and of pattern 231 in the permutation $\pi$, respectively.
We introduce the context-free grammar on the variable set $V=\{x,y,z,w,u,v\}$:
\begin{equation}\label{grammar}
  G\colon x \rightarrow  xy,\, y \rightarrow  zu,\,
       z \rightarrow  zw,\,  w \rightarrow  xv,\,
       u \rightarrow xyz^{-1}v,\,
       v \rightarrow x^{-1}zwu,
\end{equation}
which can be
used to derive  an explicit formula for the generating function of $P_n( x, y, z, w, u, v)$ in terms of  the parabolic cylinder functions.

Let us recall some basic  knowledges of the parabolic cylinder functions, where one can refer to \cite{Whittaker, Watson, Whittaker-Watson, handbook}
for more details.
Let $n$ be a nonnegative integer, then we adapt the notation of the shifted factorial as follows:
\begin{align*}
&(a)_0=1,\\[6pt] &(a)_n=a(a+1)\cdots(a+n-1),\quad n>0.
\end{align*}
The confluent hypergeometric function $_1F_1$ with parameters $a$, $b$ and $z$ is defined as
\begin{equation*}
  _1F_1(a;b;z)=\sum_{n=0}^{\infty}\frac{(a)_n}{(b)_n}\frac{z^n}{n!}.
\end{equation*}
The parabolic cylinder function $D_a(z)$ introduced by Whittaker and Watson \cite[\S16.5]{Whittaker-Watson} is defined as
 \begin{equation}
  D_a(z)=2^{\frac{a}{2}}\sqrt{\pi}e^{-\frac{z^2}{4}}\left(\frac{1}{\Gamma\left(\frac{1-a}{2}\right)}{_1F_1}\left(-\frac{a}{2};\frac{1}{2};\frac{z^2}{2}\right)
  -\frac{\sqrt{2}z}{\Gamma\left(-\frac{a}{2}\right)}{_1F_1}\left(\frac{1-a}{2};\frac{3}{2};\frac{z^2}{2}\right)\right). \label{daz}
\end{equation}

\begin{thm}\label{thmofGen z}
Let $\delta=\sqrt{xv- zu}$ and $\hat{\delta}=\sqrt{zu-xv}$. We have
  \begin{equation*}%\label{Gen( z)}
    \sum_{n=0}^\infty P_n( x, y, z, w, u, v)\frac{t^n}{n!}=
    \frac{ z(pq( w- y)+(\hat{\delta} ps-\delta qr))e^{\frac{ w- y}{2}t+\frac{\delta^2}{4}t^2}}
    {
   (\hat{\delta} s-q  y)
  D_{\frac{zu-yw}{\delta^2}}\left(\delta t+\frac{ w- y}{\delta}\right)
  +(pw-\delta r)
 D_{\frac{xv-yw}{\hat{\delta}^2}}\left(\hat{\delta} t+\frac{ y- w}{\hat{\delta}}\right)},
  \end{equation*}
  where
  \begin{align*}
    p&=D_{\frac{zu-yw}{\delta^2}}\left(\frac{ w- y}{\delta}\right),
    \quad
    q=D_{\frac{xv-yw}{\hat{\delta}^2}}\left(\frac{ y- w}{\hat{\delta}}\right),\\[6pt]
    r&=D_{\frac{xv-yw}{\delta^2}}\left(\frac{ w- y}{\delta}\right),
    \quad
    s=D_{\frac{zu-yw}{\hat{\delta}^2}}\left(\frac{ y- w}{\hat{\delta}}\right).
  \end{align*}
\end{thm}

It can be seen   that Theorem \ref{thmofGen z} is a refinement of Theorem \ref{futhm}.
By specializing the variables in Theorem \ref{thmofGen z},
we obtain the distribution of the total number of consecutive patterns $231$ and $321$ in permutations on $[n]$.
For $n\geq0$, let $L(n,k)$ be the number of permutations on $[n]$ with $k$ consecutive patterns $231$ and $321$,
and let
\begin{equation*}
  L_n(x)=\sum_{k\geq0}L(n,k)x^k.
\end{equation*}

\begin{thm}\label{coroofLnx}
  We have
  \begin{equation}\label{genLn(x)}
\sum_{n=0}^{\infty}L_n(x)\frac{t^n}{n!}=\frac{e^{\frac{t (t+2) (1-x)}{2}}}{1+xe^{\frac{x-1}{2}}\int_{t+1}^{1}e^{\frac{1-x}{2}s^2}ds}.
  \end{equation}
\end{thm}

As another specialization of Theorem \ref{thmofGen z},
 we obtain the joint distribution of exterior peaks of pattern 132 and of pattern 231 over permutations on $[n]$.
For $n\geq0$, denote by $T(n,i,j)$  the number of permutations on $[n]$ with $i$ exterior peaks of pattern 132 and $j$ exterior peaks of pattern 231, and let
\begin{equation*}
T_n(x,y)=\sum_{i,j \geq 0} T(n,i,j)x^i y^j.
\end{equation*}

\begin{thm}\label{corofTn}
 We have
\begin{equation}\label{genofTnxy}
\sum_{n=0}^\infty T_n (x,y) \frac{t^n}{n!}=\frac{e^{\frac{x-y}{2} t^2 }}{{_1F_1}\left(\frac{1-y}{2 (x-y)};\frac{1}{2};\frac{x-y}{2} t^2\right)- t\,{_1F_1}\left(\frac{1}{2}+\frac{1-y}{2(x-y)};\frac{3}{2};\frac{x-y}{2} t^2\right) }.
\end{equation}
\end{thm}

Specializations of Theorem \ref{corofTn} yield
the generating function of the number of permutations on $[n]$ with no exterior peak of pattern $132$
and the generating function of the number of permutations on $[n]$ avoiding the consecutive pattern $231$
due to  Kitaev \cite{Kitaev-2003,Kitaev-2005},
as well as the distribution of consecutive patterns $231$ in permutations on $[n]$ given by Elizalde and Noy \cite{Elizalde-Noy 2003}.

Meanwhile, by using the context-free grammar $G$,
we can derive the joint distribution of more permutation statistics.
Given a permutation $\pi=\pi_1\pi_2\cdots\pi_n$ on $[n]$,
first set $\pi_0=\pi_{n+1}=0$.
Following the terminology in \cite[\S 1.5]{Stanley-2012}, for $1\leq i\leq n$,
we call
an index $i$ a peak (or a maxima  \cite{Carlitz-1974}, or a modified maximum  \cite{GJ-2004}) if $\pi_{i-1}<\pi_i>\pi_{i+1}$,
a valley if $\pi_{i-1}>\pi_i<\pi_{i+1}$,
a double rise if $\pi_{i-1}<\pi_i<\pi_{i+1}$,
or a double descent if $\pi_{i-1}>\pi_i>\pi_{i+1}$.
For example, the permutation $4356721$ has
two peaks $1$ and $5$, one valley $2$, two double rises $3$ and $4$,
and two double descents $6$ and $7$.
Denote by $p(\pi)$, $v(\pi)$, $dd(\pi)$ and $dr(\pi)$ the numbers of peaks, valleys, double descents and double rises in permutation $\pi$, respectively.
Carlitz and Scoville \cite{Carlitz-1974} studied   the joint distribution of peaks, valleys, double rises and double descents over permutations on $[n]$.
For $n\geq1$, define
\begin{equation}\label{CSFn}
  F_n(x,y,z,w)=\sum_{\pi\in\mathfrak{S}_n}x^{p(\pi)-1}z^{v(\pi)}y^{dd(\pi)}w^{dr(\pi)}.
\end{equation}
They obtained the following generating function of $F_n(x,y,z,w)$,
see also \cite[Exercise 3.3.46]{GJ-2004} and \cite[Exercise 1.61]{Stanley-2012}.

\begin{thm}\emph{\bf(Carlitz-Scoville \cite{Carlitz-1974})}\label{CSthm}
  We have
  \begin{equation*}
    \sum_{n=1}^\infty F_n(x,y,z,w)\frac{t^n}{n!}=\frac{e^{\beta t}-e^{\alpha t}}{\beta e^{\alpha t}-\alpha e^{\beta t}},
  \end{equation*}
  where $\alpha\beta=xz$ and $\alpha+\beta=y+w$.
\end{thm}

%With some grammatical computations on the generating function
%\[Q_n(x,y,z,w)=\sum_{\pi\in\mathfrak{S}_n}x^{p(\pi)}y^{dd(\pi)}z^{p(\pi)}w^{n+1-2p(\pi)-dd(\pi)},\]
%the first named author \cite{Fu-2018} reproved Theorem \ref{CSthm}.

We also obtain a refinement of the generating function $F_n(x,y,z,w)$ by considering the patterns of peaks.
Given a permutation $\pi=\pi_1\pi_2\cdots\pi_n$ and $1\leq i\leq n$,
if $i$ is a peak with $\pi_{i-1}<\pi_i>\pi_{i+1}$,
then we call $i$ a
peak of pattern 132 if $\pi_{i-1}\leq\pi_{i+1}$,
or a peak of pattern 231  if $\pi_{i-1}>\pi_{i+1}$.
Note that only when $\pi=1$ a peak $i$ can be of pattern 132 with $\pi_{i-1}=\pi_{i+1}$.
For $n\geq1$, let
$Q_n(i,j,k,\ell)$ be the number of permutations on $[n]$ with
$i$ peaks of pattern 132, $j$ peaks of pattern 231,
$k$ double descents and $\ell$ double rises.
For $n\geq1$, define
\begin{equation*}
  Q_n( x, y, z, w, u, v)=\sum_{i,j,k,\ell=0}^n
Q_n(i,j,k,\ell) x^i  v^i  u^j  z^j y^k  w^\ell,
\end{equation*}
or equivalently,
\begin{equation}\label{qnijkl}
  Q_n( x, y, z, w, u, v)=\sum_{\pi\in\mathfrak{S}_n}
  x^{p_1(\pi)}  v^{p_1(\pi)} u^{p_2(\pi)} z^{p_2(\pi)}  y^{dd(\pi)}  w^{dr(\pi)},
\end{equation}
where $p_1(\pi)$ is the number of peaks of pattern 132 in permutation $\pi$
and $p_2(\pi)$ is the number of peaks of pattern 231
in permutation $\pi$.
Based on the same context-free grammar $G$ \eqref{grammar},
we derive the generating function of $Q_n( x, y, z, w, u, v)$.

\begin{thm}\label{thmofGen v}
   Set $Q_0(x,y,z,w,u,v)=w$. Let $\delta,\hat{\delta}, p,q,r,s$ be defined as in Theorem \ref{thmofGen z}.
   Then
   \begin{align*}
   \sum_{n=0}^\infty& Q_n( x, y, z, w, u, v)\frac{t^n}{n!} \\[6pt]
   &\qquad\quad=\frac{(\delta^2t+ w- y)(p w-\delta r)D_{\frac{xv-yw}{\hat{\delta}^2}}\left(\hat{\delta} t+\frac{ y- w}{\hat{\delta}}\right)}
     {(\hat{\delta} s-q  y)
   D_{\frac{zu-yw}{\delta^2}}\left(\delta t+\frac{ w- y}{\delta}\right)
   +(pw-\delta r)
   D_{\frac{xv-yw}{\hat{\delta}^2}}\left(\hat{\delta} t+\frac{ y- w}{\hat{\delta}}\right)}\\[6pt]
     &\qquad\qquad\qquad\quad+\frac{\delta
   (\hat{\delta} s-q  y)
   D_{\frac{xv-yw}{\delta^2}}\left(\delta t+\frac{ w- y}{\delta}\right)
   +\hat{\delta}(pw-\delta r)
   D_{\frac{zu-yw}{\hat{\delta}^2}}\left(\hat{\delta} t+\frac{ y- w}{\hat{\delta}}\right)}
     {(\hat{\delta} s-q  y)
   D_{\frac{zu-yw}{\delta^2}}\left(\delta t+\frac{ w- y}{\delta}\right)
   +(pw-\delta r)
   D_{\frac{xv-yw}{\hat{\delta}^2}}\left(\hat{\delta} t+\frac{ y- w}{\hat{\delta}}\right)}. \notag
   \end{align*}
\end{thm}

Note that Theorem \ref{thmofGen v} can be viewed as a refinement of Theorem \ref{CSthm}.
Combining Theorem \ref{thmofGen z} and Theorem \ref{thmofGen v}, we are led to the joint distribution
of exterior peaks of pattern 132
and of pattern 231
over alternating permutations on $[n]$.
For $n\geq0$, let $T^A(n,i,j)$ be the number of alternating  permutations on $[n]$ with
$i$ exterior peaks of pattern 132 and $j$ exterior peaks of pattern 231,
and let
$$
T^A_n(x,y)=\sum_{i,j \geq 0} T^A(n,i,j)x^i y^j.
$$

\begin{thm}\label{coroofTA}
We have
\begin{equation}\label{TA}
\sum_{n=0}^\infty T^A_n (x,y)\frac{t^n}{n!}
=\frac{e^{\frac{x-y}{2}t^2}\left(1+t\,{_1F_1}\left(\frac{x}{2(x-y)};\frac{3}{2};-\frac{x-y}{2}t^2\right)\right)}{{_1}F_{1}\left(-\frac{y}{2(x-y)};\frac{1}{2};\frac{x-y}{2}t^2\right)}.
\end{equation}
\end{thm}

The rest of this paper is organized as follows.
In Section \ref{subsec z},
we give an overview of the formal derivative with
respect to a context-free grammar, and
provide   grammatical labelings on permutations to
generate the polynomials $P_n(x, y, z, w, u, v)$ and $Q_n( x, y, z, w, u, v)$. In Section \ref{secofgen},
we give proofs of Theorem \ref{thmofGen z} and Theorem \ref{thmofGen v} by establishing a parabolic cylinder differential equation.
Then we show that Theorem \ref{futhm} and Theorem \ref{CSthm} are specializations of Theorem \ref{thmofGen z} and Theorem \ref{thmofGen v}.
Section \ref{secofspe} is devoted to the proofs of
Theorem \ref{coroofLnx}, Theorem \ref{corofTn} and Theorem \ref{coroofTA}.

\section{Grammatical labelings}\label{subsec z}

In this section, we first recall some basic backgrounds of
the formal derivative with respect a context-free grammar
and the grammatical labeling.
Then with the context-free grammar $G$ \eqref{grammar},
we give the corresponding grammatical labeling to
derive the generating function $P_n(x,y,z,w,u,v)$.
By using a different grammatical labeling based on the same
grammar $G$ \eqref{grammar}, we derive the generating function $Q_n(x,y,z,w,u,v)$.

Let $V$ be a variable set. A context-free grammar $G$ is defined as a set of substitution rules replacing a variable in $V$ by a Laurent polynomial of variables in $V$.
For variables $u,v\in V$, define a linear operator $D$ with respect to $G$ with the following properties:
\begin{enumerate}[(i)]
\item $D(u+v)=D(u)+D(v)$;
\item $D(uv)=D(u)v+uD(v)$;
\item $D(c)=0$, if $c$ is a constant.
\end{enumerate}
Thus, for integer $n\geq0$, the operator $D$ satisfies the Leibniz rule
 \begin{equation}\label{leibniz}
 D^n(uv)=\sum_{k=0}^n\binom{n}{k} D^k(u)D^{n-k}(v).
 \end{equation}
For a Laurent polynomial $w$ of variables in $V$, we define the generating function of  $w$ as
\begin{equation}\label{defofgen}
\gen(w,t)=\sum_{n=0}^\infty D^n(w)\frac{t^n}{n!}.
\end{equation}
Then due to \eqref{leibniz} and \eqref{defofgen}, the following relations hold:
\begin{align}
 \gen(uv,t) &= \gen(u,t)\gen(v,t), \label{product}\\[6pt]
 \gen'(u,t) &= \gen(D(u), t), \label{derivative}
 \end{align}
where $u$, $v$ are Laurent polynomials of variables in $V$, and $\gen'(u,t)$ is the normal derivative with respect to $t$.
We call the above operator $D$ the formal derivative with respect to the grammar $G$.

The idea of using the formal derivative with respect to
 a context-free grammar to study
combinatorial structures was initiated by Chen \cite{Chen-1993}.
Dumont \cite{Dumont-1996} later found grammars
for several classical combinatorial structures.
For example,  Dumont defined a context-free grammar
\begin{equation}\label{g1}
x\rightarrow xy,\, y\rightarrow xy
\end{equation}
and used the corresponding formal derivative
to generate the  Eulerian polynomials $A_n(x)$, namely, $D^n(x)|_{y=1}=xA_n(x)$.

The concept of the grammatical labeling was introduced
in \cite{Chen&Fu-2017} to build the connections between context-free grammars and the combinatorial structures.
A grammatical labeling is
 an assignment of the underlying elements of
a combinatorial structure with constants or variables, which is consistent with
 the substitution rules of a context-free grammar $G$.
For example, by using the grammar
\begin{equation}\label{g2}
x\rightarrow xy, y \rightarrow x^2,
\end{equation}
we may label the elements of a permutation $\pi$ on $[n]$ by
assigning $\pi_i$ and $\pi_{i+1}$ both label $x$ if $i$ is an exterior peak
and assigning all other elements  label $y$.
By this grammatical labeling,
one can prove that $D^n(x)$ is the generating function of
the number of permutations on $[n]$ with a given number of  exterior peaks.
This example can be found in \cite{Chen&Fu-2017},
and see \cite{Ma-2012} for a similar example but related to the normal derivative given by Ma.
In \cite{Fu-2018},
by the corresponding grammatical labeling consistent with the grammar \eqref{g3}:
\begin{equation*}
x\rightarrow xy,\,y\rightarrow xz,\,
z\rightarrow zw,\,w\rightarrow xz
\end{equation*}
the first named author gave the generating function of  $P_n(x,y,z,w)$.

%then showed
%\begin{equation*}
%  \sum_{n\geq0}^\infty P_n(x,y,z,w)\frac{t^n}{n!}=\gen(z,t),
%\end{equation*}
%where $\gen(z,t)$ is corresponding generating function to the grammar \eqref{g3}.
%Given a permutation $\pi=\pi_1\cdots\pi_n$ on $[n]$, set $\pi_0=0$,
%then for $1\leq i\leq n-1$,
%we  call $i$ an \emph{exterior peak of pattern 132} if $\pi_{i-1}<\pi_i>\pi_{i+1}$ and $\pi_{i-1}<\pi_{i+1}$,
%an \emph{exterior peak of pattern 231} if $\pi_{i-1}<\pi_i>\pi_{i+1}$ and  $\pi_{i-1}>\pi_{i+1}$,
%or a \emph{proper double descent} if $\pi_{i-1}>\pi_i>\pi_{i+1}$.
%The distribution of the proper double descent was
%studied by Elizalde-Noy \cite{Elizalde-Noy 2003}, Barry \cite{Barry-2014} and Basset \cite{Basset-2014}.
%Denote by $P_n(i,j,k)$ the number of permutations on $[n]$
%with $i$  exterior peaks of pattern 132,
%$j$ exterior peaks of pattern 231,
%and $k$  proper double descents.
% We can state the generating function of $P_n(i,j,k)$ as follows

Recall that $P_n( x, y, z, w, u, v)$ is defined in \eqref{pnijk}, that is,
\[
P_n( x, y, z, w, u, v)=\sum_{i,j,k=0}^{n}P_n(i,j,k) x^i v^i u^j z^{j+1} y^k w^{n-2i-2j-k}.
\]
Let $G$ be the context-free grammar over the variable set $V=\{ x, y, z, w, u, v\}$ defined as in \eqref{grammar}:
\begin{equation*}
  G\colon x \rightarrow  xy,\, y \rightarrow  zu,\,
       z \rightarrow  zw,\,  w \rightarrow  xv,\,
       u \rightarrow xyz^{-1}v,\,
       v \rightarrow x^{-1}zwu.
\end{equation*}
Notice that the  grammar $G$ is a refinement of grammars \eqref{g1}, \eqref{g2} and \eqref{g3}.
More precisely,
substituting $w,u$ by $x$ and $z,v$ by $y$ in the grammar $G$ reduces to the grammar \eqref{g1};
substituting $z,u,v$ by $x$ and $w$ by $y$ in the grammar $G$ reduces to the grammar \eqref{g2};
and substituting $v$ by $z$ and $u$ by $x$  in the grammar $G$ reduces to the grammar \eqref{g3}.

\begin{thm}\label{thmofP} Let $D$ be the formal derivative with respect to the grammar $G$ \eqref{grammar}.
  For $n\geq0$, we have
  \begin{equation*}
  D^n(z)=P_n( x, y, z, w, u, v).
  \end{equation*}
\end{thm}

The theorem  asserts that
the generating function $P_n( x, y, z, w, u, v)$ can be grammatically acquired
by computing the formal derivatives of $z$ with respect to the grammar $G$ \eqref{grammar}.  For instance,
\begin{align*}
%  D( z)&= zw, \\[6pt]
%  D^2( z)&= xzv+ zw^2, \\[6pt]
%  D^3( z)&=3xzwv+ z^2wu+ xyzv+ zw^3,\\[6pt]
  D^4( z)=6xzw^2v+5z^2w^2u+5xyzwv+ yz^2wu+ xy^2zv+3x^2zv^2+2xz^2uv+ zw^4.
\end{align*}
Since the coefficient of $xyzwv$
in $D^4(z)$ is $5$,
we can deduce that there are five permutations on $\{1,2,3,4\}$
with one exterior peak of pattern 132
and one proper double descent.
It is easy to check that they are 2431, 1421, 4213, 4312 and 3214.

To prove the theorem, given a permutation $\pi=\pi_1\pi_2\cdots\pi_n$ on $[n]$, we assign a labeling of $\pi$ as follows.
We first add an element $0$ at the end of $\pi$ and label it by $z$. Then for $1\leq i\leq n-1$, if $i$ is an exterior peak of pattern 132,
label $\pi_i$ by $x$ and $\pi_{i+1}$ by $v$,
if $i$ is an exterior peak of pattern 231,
label $\pi_i$ by $u$ and $\pi_{i+1}$ by $z$,
if $i$ is a proper double descent,
label $\pi_{i+1}$ by $y$.
The other elements of $\pi$ are all assigned label $w$. Denote by $w(\pi)$ the weight of permutation $\pi$
as the product of all labels of elements of  $\pi$, i.e.,
\begin{equation}\label{weight}
  w(\pi)=x^{{ep}_1(\pi)}v^{{ep}_1(\pi)} u^{{ep}_2(\pi)} z^{{ep}_2(\pi)+1} y^{{pdd}(\pi)} w^{n-2(ep_1(\pi)+ep_2(\pi))-pdd(\pi)}.
\end{equation}

\begin{example}
Let  $\pi=534621$. The labeling of $\pi$ is as follows:
\begin{equation}\label{labelPexample}
\begin{array}{ccccccc}
  5 & 3 & 4 & 6 & 2 & 1 & 0\\[-2pt]
  x & v & w & u & z & y & z
\end{array}
\end{equation}
and $w(\pi)= x y z^2 w u v$.
\end{example}

\noindent{\it Proof of Theorem \ref{thmofP}.}
We proceed to prove by induction on $n$.
By \eqref{pnijk2} and \eqref{weight},
we observe that
\begin{equation*}
  P_n( x, y, z, w, u, v)
  =\sum_{\pi\in\mathfrak{S}_n}w(\pi).
\end{equation*}
For $n=0$, the grammatical labeling of the empty permutation is given by
\begin{equation*}
  \begin{array}{c}
  0 \\[-2pt]
  z
\end{array},
\end{equation*}
which leads to $P_0( x, y, z, w, u, v)= z$.
From the view of the formal derivative,
we see that $D^0( z)= z$.
Thus the theorem holds when $n=0$.
Suppose that the theorem is valid for $n\geq1$, that is,
$D^n(z)=\sum_{\pi\in\mathfrak{S}_n}w(\pi)$.

Let $\pi=\pi_1\pi_2\ldots\pi_n$ be a permutation on $[n]$
with $i$ exterior peaks of pattern 132,
$j$ exterior peaks of pattern 231,
and $k$ proper double descents.
Hence the weight of $\pi$ is
$ x^i v^i u^j z^{j+1} y^k w^{n-2i-2j-k}$.
To complete the proof for the case  $n+1$,
we first add an element $0$ at the end of $\pi$,
then generate permutations on $[n+1]$ by
inserting $n+1$ before each element of $\pi$ (including the element 0).
Depending on the position $n+1$ inserted,
we can label $n+1$ by $x$, $u$ or $w$,
and adjust labels of some related elements.
There are six cases as follows.

\noindent{\bf Case 1.}
Insert $n+1$ before some element $\pi_\ell$ labeled by $x$.
As shown below, the position of $n+1$ becomes an exterior peak of pattern 132
since $\pi_{\ell-1}<\pi_\ell$,
so that we label $n+1$ by $x$ and relabel $\pi_\ell$ by $v$.
The position of $\pi_\ell$ becomes a proper double descent since $n+1>\pi_\ell>\pi_{\ell+1}$,
so we relabel $\pi_{\ell+1}$ by $y$.
\[\begin{array}{ccccc}
    \pi_{\ell-1}&<&\pi_\ell&>&\pi_{\ell+1}\\[-2pt]
             & & x  & & v
    \end{array}
   \quad
   \Longrightarrow
   \quad
    \begin{array}{ccccccc}
    \pi_{\ell-1}&<&n+1&>&\pi_\ell    &>&\pi_{\ell+1}\\[-2pt]
             & & x& & v& & y
    \end{array}.\]
For example, inserting $7$ before the element 5 in \eqref{labelPexample} yields the grammatical labeling:
\[
\begin{array}{cccccccc}
  7 & 5 & 3 & 4 & 6 & 2 & 1 & 0\\[-2pt]
  x & v & y & w & u & z & y & z
\end{array}.
\]

Therefore, this insertion of $n+1$ corresponds to
the substitution rule $ x\rightarrow  x y$ in grammar $G$ \eqref{grammar}.
Since there are $i$ exterior peaks of pattern 132 in $\pi$,
the insertion produces $i$ permutation on $[n+1]$
with $i$ exterior peaks of pattern 132,
$j$ exterior peaks of pattern 231,
and $k+1$ proper double descents.
The total weight of these $i$ permutations is
\[i x^i v^i u^j z^{j+1} y^{k+1} w^{n-2i-2j-k}.\]

\noindent{\bf Case 2.}
Insert $n+1$ before some element $\pi_\ell$ labeled by $v$.
We label $n+1$ by $u$ and relabel $\pi_\ell$ by $z$
since the position of $n+1$ is an exterior peak of pattern 231
due to $\pi_{\ell-1}>\pi_\ell$,
in addition, we relabel $\pi_{\ell-1}$ by $w$ since $\pi_{\ell-1}<n+1$, as illustrated below:
\[\begin{array}{ccccc}
    \pi_{\ell-2}&<&\pi_{\ell-1}&>&\pi_\ell\\[-2pt]
             & & x      & & v
    \end{array}
   \quad
   \Longrightarrow
   \quad
    \begin{array}{ccccccc}
    \pi_{\ell-2}&<&\pi_{\ell-1}&<&n+1      &>&\pi_\ell\\[-2pt]
             & & w      & & u& & z
    \end{array}.\]
For example, if we insert 7 before element 3 in \eqref{labelPexample},
then the labeling is changed to
\[
\begin{array}{cccccccc}
  5 & 7 & 3 & 4 & 6 & 2 & 1 & 0\\[-2pt]
  w & u & z & w & u & z & y & z
\end{array}.
\]

Note that this insertion corresponds to the rule $ v \rightarrow x^{-1}zwu$,
and generates $i$ permutations on $[n+1]$ with
$i-1$ exterior peaks of pattern 132,
$j+1$ exterior peaks of pattern 231,
and $k+1$ proper double descent.
The total weight of these $i$ permutations is
\[i x^{i-1} v^{i-1} u^{j+1} z^{j+2} y^{k+1} w^{n-2i-2j-k}.\]

\noindent{\bf Case 3.}
Insert $n+1$ before some element $\pi_\ell$ with label $u$.
Then the position of $n+1$ becomes
an exterior peak of pattern 132
and the position of $\pi_\ell$ becomes a proper double descent,
so that we label $n+1$ by $x$,
and relabel $\pi_\ell$ by $v$ and $\pi_{\ell+1}$ by $y$:
\[\begin{array}{ccccc}
    \pi_{\ell-1}&<&\pi_{\ell}  &>&\pi_{\ell+1}\\[-2pt]
             & & u& & z
    \end{array}
   \quad
   \Longrightarrow
   \quad
    \begin{array}{ccccccc}
    \pi_{\ell-1}&<&n+1&>&\pi_\ell    &>&\pi_{\ell+1}\\[-2pt]
             & & x& & v& & y
    \end{array}.\]
For example, we insert 7 before the element 6 in \eqref{labelPexample}. Then we obtain
\[
\begin{array}{cccccccc}
  5 & 3 & 4 & 7 & 6 & 2 & 1 & 0\\[-2pt]
  x & v & w & x & v & y & y & z
\end{array}.
\]

Hence the insertion is consistent with the rule $  u \rightarrow xz^{-1}yv$,
which produces $j$ permutations on $[n+1]$ with $i+1$ exterior peaks of pattern 132,
$j-1$ exterior peaks of pattern 231
and $k+1$ proper double descents.
The total weight of these $j$ permutations is \[j x^{i+1} v^{i+1} u^{j-1} z^{j} y^{k+1} w^{n-2i-2j-k}.\]

\noindent{\bf Case4.}
Insert $n+1$ before some element labeled by $z$.
There exist two subcases.
If $n+1$ is inserted before 0, then we  label $n+1$ by $w$:
\[\begin{array}{cc}
    \pi_n & 0 \\[-2pt]
          &  z
    \end{array}
   \quad
   \Longrightarrow
   \quad
    \begin{array}{cccc}
    \pi_{n}&<&n+1& 0\\[-2pt]
           & & w&  z
    \end{array}.\]
For example, inserting 7 before element 0 in \eqref{labelPexample}, we get
\[
\begin{array}{cccccccc}
  5 & 3 & 4 & 6 & 2 & 1 & 7 & 0\\[-2pt]
  x & v & w & u & z & y & w & z
\end{array}.
\]

If $n+1$ is inserted before $\pi_\ell$ with label $z$ for $1\leq \ell \leq n-1$,
then the label of $n+1$ is assigned $u$,
and the labele of $\pi_{\ell-1}$ is adjusted to $w$,
as showed below:
\[\begin{array}{ccccc}
    \pi_{\ell-2}&<&\pi_{\ell-1}  &>&\pi_{\ell}\\[-2pt]
             & & u& & z
    \end{array}
   \quad
   \Longrightarrow
   \quad
    \begin{array}{ccccccc}
    \pi_{\ell-2}&<&\pi_{\ell-1}&<&n+1      &>&\pi_\ell \\[-2pt]
             & & w      & & u& & z
    \end{array}.\]
For example, inserting 7 before the element 2 in \eqref{labelPexample} yields
\begin{equation*}
\begin{array}{cccccccc}
  5 & 3 & 4 & 6 & 7 & 2 & 1 & 0\\[-2pt]
  x & v & w & w & u & z & y & z
\end{array}.
\end{equation*}

In conclusion, the insertions in this case always correspond
to the substitution rule $ z\rightarrow  z w$
and yields $j+1$ permutations on $[n+1]$
with $i$ exterior peaks of pattern 132,
$j$ exterior peaks of pattern 231,
and $k$ proper double descents.
Thus, we obtain the total weight
\[(j+1) x^i v^i u^j z^{j+1} y^k w^{n-2i-2j-k+1}.\]

\noindent{\bf Case 5.}
Insert $n+1$ before some element $\pi_\ell$ with label $w$.
By the labeling rule, we know that $\pi_{\ell-1}<\pi_\ell$.
Thus the position of $n+1$ becomes an exterior peak of pattern 132
then we label $n+1$ by $x$ and relabel $\pi_\ell$ by $v$:
\[\begin{array}{ccc}
    \pi_{\ell-1}&<&\pi_{\ell}\\[-2pt]
             & & w
    \end{array}
   \quad
   \Longrightarrow
   \quad
    \begin{array}{ccccc}
    \pi_{\ell-1}&<&n+1 &>&\pi_\ell     \\[-2pt]
             & & x & & v
    \end{array}.\]
For example, we insert 7 before the element 4 in \eqref{labelPexample}, which changes the grammatical labeling to
\[
\begin{array}{cccccccc}
  5 & 3 & 7 & 4 & 6 & 2 & 1 & 0\\[-2pt]
  x & v & x & v & u & z & y & z
\end{array}.
\]

This insertion corresponds to the rule $ w\rightarrow xv$,
which generates $n-2i-2j-k$ permutations on $[n+1]$
with $i+1$ exterior peaks of pattern 132,
$j$ exterior peaks of pattern 231
and $k$ proper double descents.
Thus the total weight of these permutations equals
\[(n-2i-2j-k) x^{i+1} v^{i+1} u^j z^{j+1} y^k w^{n-2i-2j-k-1}.\]

\noindent{\bf Case 6.}
Insert $n+1$ before some element $\pi_\ell$ labeled by  $y$.
It is easy to check that
the position of $n+1$ becomes an exterior peak of pattern 231 since $\pi_{\ell-1}>\pi_\ell$.
Thus we label $n+1$ by $u$ and relabel $\pi_\ell$ by $z$:
\[\begin{array}{ccccc}
    \pi_{\ell-2}&>&\pi_{\ell-1}&>&\pi_{\ell}\\[-2pt]
             & &         & & y
    \end{array}
   \quad
   \Longrightarrow
   \quad
    \begin{array}{ccccccc}
    \pi_{\ell-2}&>&\pi_{\ell-1}&<&n+1      &>&\pi_\ell \\[-2pt]
             & &         & & u& & z
    \end{array}.\]
For example, inserting 7 before the element 1 in \eqref{labelPexample} gives
\[
\begin{array}{cccccccc}
  5 & 3 & 4 & 6 & 2 & 7 & 1 & 0\\[-2pt]
  x & v & w & u & z & u & z & z
\end{array}.
\]
Hence, this insertion is coincident with the grammatical rule $ y\rightarrow z u$
and produces $k$ permutations on $[n+1]$ with
$i$ exterior peaks of pattern 132,
$j+1$ exterior peaks of pattern 231
and $k-1$ proper double descents.
The total weight of these permutations amounts to
\[k x^i v^i u^{j+1} z^{j+2} y^{k-1} w^{n-2i-2j-k}.\]

Taking all the above six cases into account, we deduce that
\begin{align*}
  &P_{n+1}( x, y, z, w, u, v)=\sum_{\pi\in\mathfrak{S}_{n+1}}w(\pi)\\[6pt]
  & \quad = \sum_{i,j,k=0}^nP_n(i,j,k)\left(i x^i v^i u^j z^{j+1} y^{k+1} w^{n-2i-2j-k}
  +i x^{i-1} v^{i-1} u^{j+1} z^{j+2} y^{k+1} w^{n-2i-2j-k}\right.\\[6pt]
  &\quad \quad \quad\left.+j x^{i+1} v^{i+1} u^{j-1} z^{j} y^{k+1} w^{n-2i-2j-k}
  +(j+1) x^i v^i u^j z^{j+1} y^k w^{n-2i-2j-k+1}\right.\\[6pt]
  &\qquad\quad\quad\left.+(n-2i-2j-k) x^{i+1} v^{i+1} u^j z^{j+1} y^k w^{n-2i-2j-k-1}+k x^i v^i u^{j+1} z^{j+2} y^{k-1} w^{n-2i-2j-k}\right).
\end{align*}
By the grammar \eqref{grammar}, we have that
\begin{align*}
  &D( x^i v^i u^j z^{j+1} y^k w^{n-2i-2j-k})\\[6pt]
  &\quad =i x^i v^i u^j z^{j+1} y^{k+1} w^{n-2i-2j-k}
  +i x^{i-1} v^{i-1} u^{j+1} z^{j+2} y^{k+1} w^{n-2i-2j-k}\\[6pt]
  &\quad \quad \quad+j x^{i+1} v^{i+1} u^{j-1} z^{j} y^{k+1} w^{n-2i-2j-k}
  +(j+1) x^i v^i u^j z^{j+1} y^k w^{n-2i-2j-k+1}\\[6pt]
  &\qquad \quad \quad+(n-2i-2j-k) x^{i+1} v^{i+1} u^j z^{j+1} y^k w^{n-2i-2j-k-1}
  +k x^i v^i u^{j+1} z^{j+2} y^{k-1} w^{n-2i-2j-k}.
\end{align*}
Therefore,
\begin{equation*} %\label{pdn}
  P_{n+1}( x, y, z, w, u, v)=\sum_{\pi\in\mathfrak{S}_{n}}D(w(\pi)).
\end{equation*}
On the other hand, it follows from induction hypothesis  that
\[
  D^{n+1}( z)=D(D^n(z))=D\left(\sum_{\pi\in\mathfrak{S}_{n}}w(\pi)\right)
  =\sum_{\pi\in\mathfrak{S}_{n}}D(w(\pi)).
\]
Thus the theorem  holds for $n+1$, which completes the proof.
\qed

By considering a different grammatical labeling on permutations,
which is also consistent with the substitution rules of the grammar $G$ \eqref{grammar},
we give a grammatical interpretation for the generating function $Q_n(x,y,z,w,u,v)$.

\begin{thm}\label{thmofQ}
Let $D$ be the formal derivative with respect to the grammar $G$ \eqref{grammar}. For $n\geq1$, we have
\begin{equation*}
  D^n( w)
=Q_n( x, y, z, w, u, v).
\end{equation*}
\end{thm}

\noindent{\emph{Proof of Theorem \ref{thmofQ}}.}
The proof is similar to that of  Theorem \ref{thmofP},
so we only present the corresponding grammatical labeling.
Given a permutation $\pi=\pi_1\pi_2\cdots\pi_n$ on $[n]$,
first add an element $\pi_{n+1}=0$ at the end of $\pi$.
Then for $1\leq i\leq n$, if $i$ is a peak of pattern 132, we label $\pi_i$ by $x$ and $\pi_{i+1}$ by $v$,
if $i$ is a peak of pattern 231, we label $\pi_{i}$ by $u$ and $\pi_{i+1}$ by $z$,
if $i$ is a double descent, we label $\pi_{i+1}$ by $y$,
finally, if $i$ is a double rise, we label $\pi_i$ by $w$.

We claim that
this grammatical labeling is well-defined, that is,
for $1\leq i\leq n+1$, each element $\pi_i$ is assigned
exactly one label.
Since $n=1$ is trivial,  we assume $n\geq2$.
If $i=1$, then $\pi_1$ is either a peak of pattern 132
or a double rise since $0=\pi_0<\pi_1$.
Thus $\pi_1$ is labeled by $x$ or $w$.
If $i=n+1$, since $\pi_{n}$ is either a peak of pattern 231 or a double descent due to $\pi_n>\pi_{n+1}=0$,
we may label $\pi_{n+1}$ by either $z$ or $y$.
Suppose $2\leq i\leq n$. If $\pi_{i-1}<\pi_i<\pi_{i+1}$,
then $\pi_i$ is  a double rise receiving label $w$,
if $\pi_{i-1}<\pi_i>\pi_{i+1}$,
then $\pi_i$ is a peak receiving  label $x$ or label $u$.
On the other hand, if $\pi_{i-2}>\pi_{i-1}>\pi_i$,
then $\pi_{i-1}$ is a double descent
so that we  label $\pi_i$ by $y$,
if $\pi_{i-2}<\pi_{i-1}>\pi_i$,
then $\pi_{i-1}$ is a peak so that
we label $\pi_i$ by either $z$ or $v$.
Following the above labeling rules,
we can always assign $\pi_i$ a label for $1\leq i\leq n+1$,
which completes the proof.
\qed

Furthermore, let $W_n(i,j,k,\ell,m)$ be the number of permutations on $[n]$ with
$i$ peaks of pattern 132, $j$ peaks of pattern 231,
$k$ double descents, $\ell$ double rises and $m-1$ valleys.
Denote by
\begin{equation*}
  W_n( x, y, z, w, u)
=\sum_{i,j,k,\ell,m=0}^n
W_n(i,j,k,\ell,m) x^i  u^j  z^{m}   y^k  w^\ell,
\end{equation*}
or equivalently,
\begin{equation}\label{wnijkml}
  W_n( x, y, z, w, u)
=\sum_{\pi\in\mathfrak{S}_n}x^{p_1(\pi)} u^{p_2(\pi)} z^{v(\pi)+1} y^{dd(\pi)} w^{dr(\pi)}.
\end{equation}
By the definitions of peaks and valleys,
it is clear that the number of peaks in a permutation is one greater than that  of valleys.
If replacing $v$ by $z$ in $Q_n(x,y,z,w,u,v)$,
then the exponent of $z$ will record the occurrences of valleys, that is,
\begin{equation}\label{QtoW}
  Q_n(x,y,z,w,u,v)|_{v=z}=W_n(x,y,z,w,u).
\end{equation}
%Hence Theorem \ref{thmofQ} also gives the generating function of $W_n(i,j,k,\ell,m)$.

%Let
%$P_0( x, y, z, w, u, v)= z$.
%By Theorem \ref{thmofP} and \eqref{defofgen}, we have
%\begin{align*}
%  \sum_{n=0}^{\infty}P_n( x, y, z, w, u, v)\frac{t^n}{n!}&=\gen(z,t).
%\end{align*}

%In section 3, we focus on grammatically computing the generating functions $\gen(z,t)$ and $\gen(v,t)$.
%By solving a general parabolic cylinder differential equation,
%we first express the generating function $\gen(x^{-1/2} z^{-1/2},t)$ in terms of the parabolic cylinder functions $D_a(z)$.
%Then using the Leibniz rule, we obtain the main theorem of this paper on the generating function $\gen(z,t)$.

%
%Let $Q_0( x, y, z, w, u, v)= w$.
%By using the substitution rule $ z\rightarrow  z w$ in grammar $G$,
%we can directly obtain the explicit form of $\gen(w,t)$ from Theorem \ref{thmofGen z}.
%
%
%In this section we introduce two grammatical labelings of permutations
%with respect to the grammar $G$ \eqref{grammar}. Based on these labelings, we give proofs of Theorem \ref{thmofP} and Theorem \ref{thmofQ}.

\begin{cor}\label{thmofQ'}
Let $D$ be the formal derivative with respect to the grammar $G$ \eqref{grammar}. For $n\geq1$, we have
\begin{align*}
  D^n( w)|_{v = z}
=W_n(x, y, z, w, u).
\end{align*}
\end{cor}

By
using the Leibniz rule \eqref{leibniz} on $D^{n+1}( z)=D^n( z w)$, we have
\[ D^{n+1}( z)=\sum_{k=0}^{n} {n \choose k} D^k( z)D^{n-k}( w).\]
Hence for $n\geq1$, combining Theorem \ref{thmofP} and Theorem \ref{thmofQ},
we can establish a convolutional relation:
  \begin{equation*}
    P_{n+1}( x, y, z, w, u, v)=\sum_{k=0}^{n}{n\choose k}
    P_k( x, y, z, w, u, v)Q_{n-k}( x, y, z, w, u, v).
  \end{equation*}

\section{Connection to the parabolic cylinder functions}\label{secofgen}

In this section,
by using the grammar
\begin{equation*}
  G\colon x \rightarrow  xy,\, y \rightarrow  zu,\,
       z \rightarrow  zw,\,  w \rightarrow  xv,\,
       u \rightarrow xyz^{-1}v,\,
       v \rightarrow x^{-1}zwu
\end{equation*}
as in \eqref{grammar}, we give the proofs of Theorem \ref{thmofGen z} and Theorem \ref{thmofGen v}.
To be more specific,
combining  Theorem \ref{thmofP} and \eqref{defofgen}, we have
\begin{equation}\label{Ptogenz}
  \sum_{n=0}^\infty P_n( x, y, z, w, u, v)\frac{t^n}{n!}=\gen(z,t).
\end{equation}
Similarly, by Theorem \ref{thmofQ} and \eqref{defofgen},
along with the assumption that $Q_0(x,y,z,w,u,v)=w$, we find that
\begin{equation}\label{Qtogenw}
  \sum_{n=0}^\infty Q_n( x, y, z, w, u, v)\frac{t^n}{n!}=\gen(w,t).
\end{equation}
Therefore, to prove Theorem \ref{thmofGen z} and Theorem \ref{thmofGen v},
 we need to compute  $\gen(z,t)$ and $\gen(w,t)$.
To this end, we first give the explicit formula for
the generating function $\gen(x^{-1/2} z^{-1/2},t)$ by solving the parabolic cylinder differential equation.
Recall that the parabolic cylinder differential equation
is a second-order ordinary differential equation of the form
  \begin{equation}\label{diequation2}
 y^{''}(z)+\left(a+\frac{1}{2}-\frac{z^2}{4}\right)y(z)=0
\end{equation}
whose solution can be written as
$$
y=c_1D_a(z) +c_2D_{-a-1}(iz),
$$
where $D_a(z)$ and $D_{-a-1}(iz)$ are independent,
and $c_1$, $c_2$ are   constants  to be determined by  the initial conditions of \eqref{diequation2}.

Based on the solutions of the equation \eqref{diequation2},
it is not difficult to solve a more general differential equation of the form
\begin{equation}\label{newadd}
 y^{''}(z)-(az^2+bz+c)y(z)=0.
\end{equation}

\begin{thm}\label{theo3.1}
For the general parabolic cylinder differential equation \eqref{newadd},
all solutions can be written as
  \begin{equation*}
    c_1D_{\frac{b^2-4ac-4a^{3/2}}{8a^{3/2}}}\left(\sqrt{2}a^{1/4}z+\frac{b}{\sqrt{2}a^{3/4}}\right)+c_2D_{\frac{4ac-b^2-4a^{3/2}}{8a^{3/2}}}\left(i\left(\sqrt{2}a^{1/4}z+\frac{b}{\sqrt{2}a^{3/4}}\right)\right),
  \end{equation*}
  where $c_1$ and $c_2$ are   constants.
\end{thm}
\pf
By equation \eqref{diequation2}, we know that
these two parabolic cylinder functions are independent,
thus we only need to prove that they are solutions of \eqref{newadd}.
To this end, we first introduce  two recurrence relations related to $D_a(z)$ and  $D_a'(z)$ in \cite[pp.16]{ParabolicCylinderD}:
  \begin{align}
    D_a'(z)&= \frac{1}{2}zD_a(z)-D_{a+1}(z), \label{recD1} \\[6pt]
    D_a'(z)&=aD_{a-1}(z)-\frac{1}{2}zD_{a}(z).\label{recD2}
%    D_a'(z)&=\frac{1}{2}\left(aD_{a-1}(z)-D_{a+1}(z)\right)\label{recD3}\\[6pt]
%    D_a(z)&=\frac{1}{z}\left(aD_{a-1}(z)+D_{a+1}(z)\right)\label{recD4}
  \end{align}
By \eqref{recD1} and \eqref{recD2}, for any two constants $r$ and $s$, we can deduce that
\begin{align}
  D'_a(rz+s)&=\frac{r}{2}(rz+s)D_a(rz+s)-rD_{a+1}(rz+s),\label{D1rz+s}\\[6pt]
  D'_{a+1}(rz+s)&=r(a+1)D_a(rz+s)-\frac{r}{2}(rz+s)D_{a+1}(rz+s).\label{D1a+1rz+s}
\end{align}
Combining \eqref{D1rz+s} and \eqref{D1a+1rz+s} yields
\begin{align}
  D''_a(rz+s)=&\left(\frac{r}{2}(rz+s)D_a(rz+s)-rD_{a+1}(rz+s)\right)'\notag\\[6pt]
  &=\frac{r^2}{2}D_a(rz+s)+\frac{r}{2}(rz+s)\left(\frac{r}{2}(rz+s)D_a(rz+s)-rD_{a+1}(rz+s)\right)\notag\\[6pt]
  &\qquad\qquad-r^2\left((a+1)D_a(rz+s)-\frac{1}{2}(rz+s)D_{a+1}(rz+s)\right)\notag\\[6pt]
  &=\frac{r^2}{4}\left((rz+s)^2-4a-2\right)D_a(rz+s).\label{D2rz+S}
\end{align}
Therefore,
plugging $$D_{\frac{b^2-4ac-4a^{3/2}}{8a^{3/2}}}\left(\sqrt{2}a^{1/4}z+\frac{b}{\sqrt{2}a^{3/4}}\right) \mbox{ and } D_{\frac{4ac-b^2-4a^{3/2}}{8a^{3/2}}}\left(i\left(\sqrt{2}a^{1/4}z+\frac{b}{\sqrt{2}a^{3/4}}\right)\right)$$
into \eqref{D2rz+S}, we see that
both of them satisfy the general parabolic cylinder differential equation \eqref{newadd}.
This completes the proof.
\qed

The following theorem states that the generating function $\gen( x^{-1/2} z^{-1/2},t)$  satisfies a general parabolic cylinder differential equation.

\begin{thm}\label{theo3.2} Let
\[ f(t)=\gen(x^{-1/2} z^{-1/2},t)\]
and
  \[ \alpha=(y+w)^2-2( x v+ z u), \quad \beta= 2 ( w- y)(x v- z u), \quad \gamma=2( x  v- z u)^2.\]
  We have
  \begin{equation}\label{formulah(t)}
    f''(t)-\left(\frac{\gamma}{8}t^2+\frac{\beta}{4}t+\frac{\alpha}{4}\right)f(t)=0.
  \end{equation}
 % with initial conditions $h(0)= x^{-1/2}  z^{-1/2}$ and $h'(0)=D( x^{-1/2} z^{-1/2})$.
\end{thm}
\pf By the grammar $G$ \eqref{grammar} and the corresponding formal derivative $D$, one can check that
\begin{align}
  %&D( x^{-1/2}  z^{-1/2})=-\frac{1}{2} x^{-1/2}  z^{-1/2}(y+w),\label{dd1} \\[6pt]
  &D^2( x^{-1/2}  z^{-1/2})=\frac{1}{4} x^{-1/2}  z^{-1/2}((y+w)^2-2(x v+z u))=\frac{\alpha }{4} x^{-1/2}  z^{-1/2}. \label{dd2}
\end{align}
Moreover, since $D(x v- z u)=0$, it follows that
\begin{equation}\label{dd5}
D(\alpha)=\beta,\quad D(\beta)=\gamma,\quad D(\gamma)=4(x v-z u)D(x v-z u)=0.
\end{equation}
By \eqref{derivative}, we have
\begin{equation}\label{addnew2}
  f^{(4)}(t)=\sum_{n=0}^{\infty} D^{n+4}( x^{-1/2}  z^{-1/2}) \frac{t^n}{n!}.
\end{equation}
From the Leibniz rule \eqref{leibniz} and \eqref{dd2}, we see that
\begin{equation*}\label{re2}
D^{n+4}( x^{-1/2}  z^{-1/2})=D^{n+2}\left(\frac{\alpha}{4}  x^{-1/2}  z^{-1/2}\right)=\frac{1}{4}\sum_{k=0}^{n+2}\binom{n+2}{k}D^k(\alpha)D^{n+2-k}( x^{-1/2}  z^{-1/2}).
\end{equation*}
By \eqref{dd5}, the summation has only three nonzero terms, implying that
\begin{align}
&D^{n+4}( x^{-1/2}  z^{-1/2})\label{dd6}\\[6pt]
&\,\,\,=\frac{\alpha}{4} D^{n+2}( x^{-1/2}  z^{-1/2})+\frac{n+2}{4}\beta D^{n+1}( x^{-1/2}  z^{-1/2})+\frac{1}{4}\binom{n+2}{2}\gamma D^{n}( x^{-1/2}  z^{-1/2}). \notag
\end{align}
Putting \eqref{dd6} into \eqref{addnew2}, and using \eqref{derivative} again, we deduce that
\begin{equation}\label{ode2'}
  f^{(4)}(t)-\left(\frac{\gamma}{8}t^2+\frac{\beta}{4}t+\frac{\alpha}{4}\right)f''(t)
  -\left(\frac{\gamma}{2}t+\frac{\beta}{2}\right)f'(t)-\frac{\gamma}{4}f(t)=0.
\end{equation}

Let $$g(t)=f''(t)-\left(\frac{\gamma}{8}t^2+\frac{\beta}{4}t+\frac{\alpha}{4}\right)f(t).$$
Taking the second-order derivative with respect to $t$ on both sides gives
%\begin{equation}\label{y'(t)}
%  g'(t)=h^{(3)}(t)-\left(\frac{\gamma}{8}t^2+\frac{\beta}{4}t+\frac{\alpha}{4}\right)h'(t)
%  -\left(\frac{\gamma}{4}t+\frac{\beta}{4}\right)h(t)
%\end{equation}
%and
\begin{equation*}
  g''(t)=f^{(4)}(t)-\left(\frac{\gamma}{8}t^2+\frac{\beta}{4}t+\frac{\alpha}{4}\right)f''(t)
  -\left(\frac{\gamma}{2}t+\frac{\beta}{2}\right)f'(t)
  -\frac{\gamma}{4}f(t).
\end{equation*}
It follows from \eqref{ode2'} that
\begin{equation*}
  g''(t)=0.
\end{equation*}
Clearly, we have $g(t)=c_1t+c_2$ with the indeterminate constants $c_1$ and $c_2$, that is,
\begin{equation}\label{formulay''(t)}
  f''(t)-\left(\frac{\gamma}{8}t^2+\frac{\beta}{4}t+\frac{\alpha}{4}\right)f(t)-c_1t-c_2=0.
\end{equation}
Set $t=0$ and recall $f^{(n)}(0)=D^n(x^{-1/2}  z^{-1/2})$.
By \eqref{dd2}, we have
\begin{equation*}
  c_2=f''(0)-\frac{\alpha}{4}f(0)=D^2( x^{-1/2}  z^{-1/2})-\frac{\alpha}{4} x^{-1/2}  z^{-1/2}=0.
\end{equation*}
Hence, \eqref{formulay''(t)} can be reduced to
\begin{equation*}
  f''(t)-\left(\frac{\gamma}{8}t^2+\frac{\beta}{4}t+\frac{\alpha}{4}\right)f(t)-c_1t=0.
\end{equation*}
Then taking the derivative with respect to $t$ on both sides yields
\begin{equation}\label{formulay''(t)2}
  f^{(3)}(t)-\left(\frac{\gamma}{8}t^2+\frac{\beta}{4}t+\frac{\alpha}{4}\right)f'(t)
  -\left(\frac{\gamma}{4}t+\frac{\beta}{4}\right)f(t)-c_1=0.
\end{equation}
Setting $t=0$ in \eqref{formulay''(t)2} gives
\begin{align}
c_1&=f^{(3)}(0)-\frac{\alpha}{4}f'(0)-\frac{\beta}{4}f(0)\notag\\[6pt]
   &=D^3(x^{-1/2}  z^{-1/2})-\frac{\alpha}{4}D( x^{-1/2}  z^{-1/2})-\frac{\beta}{4} x^{-1/2}  z^{-1/2}.\label{dd3}
\end{align}
Hence, from \eqref{dd2}, we deduce that
\begin{equation}\label{dd4}
D^3(x^{-1/2}  z^{-1/2})=D\left(\frac{\alpha }{4} x^{-1/2}  z^{-1/2}\right)=\frac{\alpha}{4}D( x^{-1/2}  z^{-1/2})+\frac{\beta}{4} x^{-1/2}  z^{-1/2}.
\end{equation}
Combining \eqref{dd3} and \eqref{dd4}, we see that $c_1=0$.
Therefore, given $c_1=c_2=0$,  equation \eqref{formulay''(t)}  reduces  to \eqref{formulah(t)}. This completes the proof.
\qed

Combing Theorem \ref{theo3.1} and Theorem \ref{theo3.2},
we give the explicit form for  the generating function $\gen( x^{-1/2} z^{-1/2},t)$ in terms of the parabolic cylinder functions.
\begin{thm}\label{thmGen( z^-1/2 x^-1/2)}
   Let $\delta,\hat{\delta}, p,q,r,s$ be defined as in Theorem \ref{thmofGen z}.
   We have
  \begin{align}
\gen&\left( x^{-1/2} z^{-1/2},t\right)\label{Gen( z^-1/2 x^-1/2)}\\[6pt]
&\qquad\qquad\quad=\frac{ x^{-1/2} z^{-1/2}}{pq( w- y)+( \hat{\delta} ps-\delta qr)}\left(
   (\hat{\delta} s-q  y)
  D_{\frac{z u-y w}{\delta^2}}\left(\delta t+\frac{ w- y}{\delta}\right)\right.\notag\\[6pt]
  &\qquad\qquad\qquad\qquad\qquad\left.+(p w-\delta r)
 D_{\frac{x v-y w}{\hat{\delta}^2}}\left(\hat{\delta} t+\frac{ y- w}{\hat{\delta}}\right)\right).\notag
  \end{align}
\end{thm}
\pf From Theorem \ref{theo3.1} and Theorem \ref{theo3.2}, $\gen( x^{-1/2} z^{-1/2},t)$ can be expressed as follows
\begin{align*}
  \gen&\left( x^{-1/2} z^{-1/2},t\right)\\[6pt]
  &\quad=c_1D_{\frac{zu-yw}{xv-zu}}\left(\sqrt{xv-zu}t+\frac{ w- y}{\sqrt{xv-zu}}\right)+c_2D_{\frac{yw-xv}{xv-zu}}\left(i\left(\sqrt{xv-zu}t+\frac{ w- y}{\sqrt{xv-zu}}\right)\right)
\end{align*}
with the indeterminate constants $c_1$ and $c_2$.
By noticing that $i \sqrt{xv-zu}=\sqrt{zu-xv}$, we rewrite $\gen\left( x^{-1/2} z^{-1/2},t\right)$ as
\begin{align*}
  \gen&\left( x^{-1/2} z^{-1/2},t\right)\\[6pt]
  &\quad=c_1D_{\frac{zu-yw}{xv-zu}}\left(\sqrt{xv-zu}t+\frac{ w- y}{\sqrt{xv-zu}}\right)+c_2D_{\frac{xv-yw}{zu-xv}}\left(\sqrt{zu-xv}t+\frac{ y- w}{\sqrt{zu-xv}}\right).
\end{align*}

Let $\delta,\hat{\delta}, p,q,r,s$ be defined as  in Theorem \ref{thmofGen z}.
Since $\gen( x^{-1/2} z^{-1/2},0)= x^{-1/2} z^{-1/2}$,
we obtain
\begin{equation}\label{eq1c1c2}
  c_1p+c_2q= x^{-1/2} z^{-1/2}.
\end{equation}
Moreover,
\[ \gen'( x^{-1/2} z^{-1/2},0)
=D( x^{-1/2} z^{-1/2})=-\frac{1}{2} x^{-1/2}
 z^{-1/2}(y+w).\]
Thus, by \eqref{D1rz+s}, we see that
\begin{equation}\label{eq2c1c2}
  c_1\left(\frac{ w- y}{2}p-\delta r\right)+
  c_2\left(\frac{ y- w}{2}q-\hat{\delta} s\right)
  =-\frac{1}{2} x^{-1/2}  z^{-1/2}(y+w).
\end{equation}
Combining \eqref{eq1c1c2} and \eqref{eq2c1c2}, we arrive at
\begin{equation*}
  c_1=\frac{(\hat{\delta} s-q y)x^{-1/2} z^{-1/2}}{pq( w- y)+( \hat{\delta} ps-\delta qr)}
  \quad \mbox{and} \quad
  c_2=\frac{(pw -\delta r) x^{-1/2} z^{-1/2}}{pq( w- y)+( \hat{\delta} ps-\delta qr)}
.
\end{equation*}
This completes the proof.
\qed

Now we are ready to prove
Theorem \ref{thmofGen z} and Theorem \ref{thmofGen v}.

\noindent{\emph{Proof of Theorem \ref{thmofGen z}}.}
Using mathematical induction on $n$  with the fact $D^2( w- y)=D( x v- z u)=0$, we see that for $n\geq0$,
\begin{equation*}
    D^{n}( x^{-1}z)= x^{-1}z\sum_{k=0}^{\lfloor n/2\rfloor}\frac{n!}{2^k(n-2k)!k!}( w- y)^{n-2k}( x  v- z u)^k.
  \end{equation*}
Therefore, we have
\begin{align}
  \gen( x^{-1}z,t)&=\sum_{n=0}^\infty
   x^{-1}z\left(\sum_{k=0}^{\lfloor n/2\rfloor}\frac{n!}{2^k(n-2k)!k!}( w- y)^{n-2k}( x  v- z u)^k\right)\frac{t^n}{n!} \nonumber \\[6pt]
  &= x^{-1}z\sum_{k=0}^\infty\sum_{n=0}^\infty\frac{( w- y)^{n}( x v- z u)^k t^{n+2k}}{2^kn!k!} \nonumber\\[6pt]
%  &= x^{-1}z\sum_{n=0}^\infty( v- u)^n\frac{t^n}{n!}\cdot\sum_{k=0}^\infty\left(\frac{ x  w- yz}{2}\right)^k\frac{(t^2)^k}{k!} \nonumber\\[6pt]
  &= x^{-1}ze^{(w- y)t+\frac{( x v- z u)t^2}{2}}.\label{Gen( z x^-1)}
\end{align}
On the other hand, by \eqref{product}, we have
\begin{equation*}
  \left(\gen(z,t)\cdot\gen( x^{-1/2}z^{-1/2},t)\right)^2=\gen(x^{-1/2} z^{1/2} ,t)^2=\gen(x^{-1} z,t),
\end{equation*}
which leads to
\begin{equation}\label{toGen z}
  \gen(z,t)=\frac{\gen(  x^{-1}z, t)^{1/2}}{\gen( x^{-1/2}z^{-1/2},t)}.
\end{equation}
Plugging \eqref{Gen( z^-1/2 x^-1/2)} and \eqref{Gen( z x^-1)} into \eqref{toGen z},
we arrive at the generating function $\gen(z,t)$.
\qed

\vspace{8pt}

\noindent{\emph{Proof of Theorem \ref{thmofGen v}}.}
From the substitution rule $ z\rightarrow  z w$ in the grammar \eqref{grammar}, utilizing \eqref{product} and \eqref{derivative},
we obtain that
\begin{equation*}
  \gen'( z,t)=\gen(D( z),t)=\gen( z w,t)=\gen(z,t)\gen(w,t),
\end{equation*}
thus,
\begin{equation}\label{getgen(w,t)}
  \gen(w,t)=\frac{\gen'( z,t)}{\gen(z,t)}.
\end{equation}
Using \eqref{D1rz+s} on $\gen(z,t)$ to compute \eqref{getgen(w,t)},
we are led to Theorem \ref{thmofGen v}.
\qed

To see that Theorem \ref{thmofGen z} implies Theorem \ref{futhm},
replacing $u$ by $x$ and $v$ by $z$ in \eqref{pnijk2}
gives \eqref{fupijk}, that is,
\begin{equation*}
  P_n(x,y,z,w,u,v)\big|_{u=x \atop v=z}=P_n(x,y,z,w).
\end{equation*}
Thus by \eqref{Ptogenz}, we derive
\begin{equation*}
  \sum_{n=0}^{\infty}P_n(x,y,z,w)\frac{t^n}{n!}=\gen(z,t)\big|_{u=x \atop v=z}.
\end{equation*}
It follows from \eqref{toGen z} that
\begin{equation}\label{tofu1}
  \gen(z,t)\big|_{u=x \atop v=z}=\frac{\Big(\gen(  x^{-1}z, t)\big|_{u=x\atop v=z}\Big)^{1/2}}{\gen( x^{-1/2}z^{-1/2},t)\big|_{u=x\atop v=z}}.
\end{equation}
From \eqref{Gen( z x^-1)}, we can directly deduce that
\begin{equation}\label{tofu2}
  \gen(  x^{-1}z, t)\big|_{u=x\atop v=z}=x^{-1}ze^{(w- y)t}.
\end{equation}
With the replacement $u\rightarrow x$ and $v\rightarrow z$,
the parabolic cylinder differential equation \eqref{formulah(t)} reduces to the second order linear differential equation
\begin{equation*}
  h''(t)-((y+w)-4xz)h(t)=0,
\end{equation*}
where $h(t)=\gen( x^{-1/2}z^{-1/2},t)\big|_{u=x\atop v=z}$.
According to  the initial conditions
$$h(0)=x^{-\frac{1}{2}}z^{-\frac{1}{2}}\quad\mbox{and}\quad h'(0)=-\frac{1}{2}x^{-\frac{1}{2}}z^{-\frac{1}{2}}(y+w),$$
solving the linear differential equation gives
\begin{align}\label{tofu3}
&\gen( x^{-\frac{1}{2}}z^{-\frac{1}{2}},t)|_{u=x\atop v=z}\\[6pt]
  &\qquad\qquad\quad=\frac{x^{-\frac{1}{2}}z^{-\frac{1}{2}}}{2}
    \left(\left(1+\frac{y+w}{\sqrt{(y+w)^2-4xz}}\right)e^{-\frac{t}{2}\sqrt{(y+w)^2-4xz}}\right.\notag\\[6pt] &\qquad\qquad\qquad\qquad\qquad\quad\left.+\left(1-\frac{y+w}{\sqrt{(y+w)^2-4xz}}\right)e^{\frac{t}{2}\sqrt{(y+w)^2-4xz}}\right).\notag
\end{align}
Hence, plugging \eqref{tofu2} and \eqref{tofu3} into \eqref{tofu1},
we arrive at
\begin{align}\label{tofu4}
  \gen(z,t)&\big|_{u=x \atop v=z}\notag\\[6pt]
  &=\frac{2z\sqrt{(y+w)^2-4xz}e^{\frac{t}{2}(w-y+\sqrt{(y+w)^2-4xz})}}{y+w+\sqrt{(y+w)^2-4xz}-(y+w-\sqrt{(y+w)^2-4xz})e^{t\sqrt{(y+w)^2-4xz}}},
\end{align}
which proves Theorem \ref{futhm}.

Next we show that Theorem \ref{thmofGen v} implies Theorem  \ref{CSthm}.
Substituting $u$ by $x$ in \eqref{wnijkml} and comparing with \eqref{CSFn},
then for $n\geq1$, we deduce
\begin{equation*}
  W_n(x,y,z,w,u)|_{u=x}=xzF_n(x,y,z,w).
\end{equation*}
Together with \eqref{QtoW}, we have
\begin{equation*}
  Q_n(x,y,z,w,u)|_{u=x\atop v=z}=xzF_n(x,y,z,w).
\end{equation*}
By \eqref{Qtogenw}, we see
\begin{align}\label{toCS1}
  \sum_{n=1}^\infty F_n(x,y,z,w)\frac{t^n}{n!}=&x^{-1}z^{-1}\sum_{n=1}^\infty Q_n(x,y,z,w,u,v)\big|_{u=x\atop v=z}\frac{t^n}{n!} \notag\\[6pt]
  =&x^{-1}z^{-1}\left(\gen(w,t)\big|_{u=x\atop v=z}-w\right).
\end{align}
It follow from \eqref{getgen(w,t)} that
\begin{equation}\label{toCS2}
  \gen(w,t)\big|_{u=x\atop v=z}=\frac{\Big(\gen( z,t)\big|_{u=x\atop v=z}\Big)'}{\gen(z,t)\big|_{u=x\atop v=z}}.
\end{equation}
Let $\alpha\beta=xz$ and $\alpha+\beta=y+w$,
then by \eqref{tofu4} and \eqref{toCS2}, we derive
\begin{equation}\label{toCS3}
  \gen(w,t)\big|_{u=x\atop v=z}=\alpha\beta\frac{e^{\beta t}-e^{\alpha t}}{\beta e^{\alpha t}-\alpha e^{\beta t}}+w.
\end{equation}
Therefore, plugging \eqref{toCS3} into  \eqref{toCS1} proves Theorem \ref{CSthm}.

%\begin{remark}
%Combining Corollary \ref{thmofQ'} and Theorem \ref{thmofGen v},
%we can give the generating function of $ W_n( x, y, z, w, u)$
%by replacing $v$ with $z$ in \eqref{Gen v}.
%\end{remark}

\section{Specializations} \label{secofspe}

In this section, we discuss several applications of Theorem \ref{thmofGen z} and Theorem \ref{thmofGen v} by specializing the variables in $\{ x, y, z, w,  u,  v\}$.
Before this, we recall the definition of  {consecutive patterns} \cite{Elizalde-Noy 2003} in permutations.
Let $m\leq n$ be two positive integers, and $\pi=\pi_1\pi_2\cdots\pi_n$ be a permutation on $[n]$
and $\sigma=\sigma_1\sigma_2\cdots\sigma_m$ be a permutation on $[m]$.
Then we say $\pi$ contains a consecutive pattern $\sigma$
if there exist $m$ consecutive elements $\pi_{i+1}\pi_{i+2}\cdots\pi_{i+m}$ in $\pi$ ($0\leq i\leq n-m$)
such that
$\rho(\pi_{i+1}\pi_{i+2}\cdots\pi_{i+m})=\sigma$,
where $\rho$ is the reduction consisting in relabeling the
elements with $\{1,2,\ldots,m\}$
so that they keep the same order relationships they have in $\pi$.

Let $L(n,k)$ be the number of permutations of $[n]$ with $k$ consecutive patterns $231$ and $321$, and let
\begin{equation}\label{Ln(x)}
  L_n(x)=\sum_{k\geq0}L(n,k)x^k.
\end{equation}
In other words, $L(n,k)$ counts the number of permutations of $[n]$
whose total number of exterior peaks of pattern 231 and proper double descents is $k$.

%\begin{thm}\label{coroofLnx}
%  We have
%  \begin{equation}\label{genLn(x)}
%\sum_{n=0}^{\infty}L_n(x)\frac{t^n}{n!}=\frac{e^{\frac{t (t+2) (1-x)}{2}}}{1+xe^{\frac{x-1}{2}}\int_{t+1}^{1}e^{\frac{1-x}{2}s^2}ds}.
%  \end{equation}
%\end{thm}

\noindent{\emph{Proof of Theorem \ref{coroofLnx}}.}
Note that by \eqref{pnijk}, we have
\begin{equation}\label{PtoL}
  P_n(1,x,1,1,x,1)=L_n(x).
\end{equation}
Due to the explicit form of $D_a(z)$ in \eqref{daz} with the fact $\lim_{z\rightarrow0}\frac{1}{\Gamma(z)}=0$,
we can easily obtain
\begin{equation}\label{D0D1}
D_0(z)=e^{-\frac{z^2}{4}} \mbox{ and } D_1(z)=ze^{-\frac{z^2}{4}}.
\end{equation}
Moreover, we see that
\begin{equation*}
  D_{-1}(z)=\sqrt{\frac{\pi}{2}} e^{-\frac{z^2}{4}}\left( e^{\frac{z^2}{2}}-\frac{\sqrt{2}z}{\sqrt{\pi}}{_1F_1}\left(1;\frac{3}{2};\frac{z^2}{2}\right)\right)=\sqrt{\frac{\pi}{2}}e^{\frac{z^2}{4}}-ze^{-\frac{z^2}{4}}{_1F_1}\left(1;\frac{3}{2};\frac{z^2}{2}\right).
\end{equation*}
Notice that
\begin{equation}\label{1f1toerf}
{_1F_1}\left(1;\frac{3}{2};z^2\right)=\frac{\sqrt{\pi}}{2z}e^{z^2}\erf(z),
\end{equation}
where $\erf(x)$ is the  {error function} defined by
  \begin{equation*}
    \erf(x)=\frac{2}{\sqrt{\pi}}\int_{0}^{x}e^{-s^2}ds.
  \end{equation*}
Then setting $z\rightarrow\frac{z}{\sqrt{2}}$ in \eqref{1f1toerf}  gives
$$
  {_1F_1}\left(1;\frac{3}{2};\frac{z^2}{2}\right)=\sqrt{\frac{\pi}{2}}z^{-1}e^{\frac{z^2}{2}}\erf\left(\frac{z}{\sqrt{2}}\right),
$$
which leads to
\begin{equation}\label{D-1}
  D_{-1}(z)=\sqrt{\frac{\pi}{2}}e^{\frac{z^2}{4}}\left(1-\erf\left(\frac{z}{\sqrt{2}}\right)\right).
\end{equation}
Therefore, by \eqref{PtoL},
setting $ x,z,w,v\rightarrow1$ and $y,u\rightarrow x$ in Theorem \ref{thmofGen z},  utilizing \eqref{D0D1} and \eqref{D-1}, we have
  \begin{equation*}
    \sum_{n=0}^{\infty}\frac{L_n(x)}{n!}t^n=\frac{ \sqrt{2(x-1)} e^{\frac{t (t+2) (1-x)}{2}}}{\sqrt{2(x-1)}+\sqrt{\pi} e^{\frac{x-1}{2}} x\left(\erf\left(\sqrt{\frac{x-1}{2}}\right)-\erf\left(\frac{ (t+1) \sqrt{x-1}}{\sqrt{2}}\right)\right)}.
  \end{equation*}
After replacing the error function $\erf(x)$ by integration, we finally arrive at \eqref{genLn(x)}.
This completes the proof.
\qed

By setting $x=0$ in \eqref{Ln(x)},
we have $L_n(0)=L(n,0)$,
which is the number of permutations on $[n]$ avoiding the consecutive patterns 231 and 321.
It follows from Theorem \ref{coroofLnx} that
\begin{equation}\label{genofLn0}
  \sum_{n=0}^{\infty}\frac{L_n(0)}{n!}t^n=e^{t+\frac{t^2}{2}}.
\end{equation}
Note that \eqref{genofLn0} is also the generating function of the number of involutions on $[n]$, see \cite[Eq. 5.32]{Stanley-2012},
which implies that
the number of permutations on $[n]$ avoiding the consecutive patterns 231 and 321
equals the number of involutions on $[n]$.
A bijective proof given by Callan can be found in \cite[A000085]{OEIS}.

Let $T(n,i,j)$ be the number of permutations on $[n]$ with $i$ exterior peaks of pattern 132 and $j$ exterior peaks of pattern 231, and let
\begin{equation}\label{Tnxy}
T_n(x,y)=\sum_{i,j \geq 0} T(n,i,j)x^i y^j.
\end{equation}

\noindent{\emph{Proof of Theorem \ref{corofTn}}.}
By \eqref{pnijk}, we know that $P_n(x,1,1,1,y,1)=T_n(x,y)$.
Thus after  replacements $y,z,w,v\rightarrow1$
and $u\rightarrow y$ in Theorem \ref{thmofGen z},
we obtain the numerator as
\begin{equation}\label{numerator}
\lim_{ y,z,w,v\rightarrow1 \atop u\rightarrow y}
  z(pq( w- y)+(\hat{\delta} ps-\delta qr))e^{\frac{ w- y}{2}t+\frac{\delta^2}{4}t^2}=\pi e^{\frac{x-y}{4}t^2}\Theta,
\end{equation}
where
\[\Theta=\frac{\sqrt{y-x}}{\Gamma\left(\frac{1-x}{2(y-x)}\right)\Gamma\left(\frac{1}{2}+\frac{1-y}{2(x-y)}\right)}-\frac{\sqrt{x-y}}{\Gamma\left(\frac{1-y}{2(x-y)}\right)\Gamma\left(\frac{1}{2}+\frac{1-x}{2(y-x)}\right)}.\]
For the denominator, we have
\begin{align*}
&\lim_{ y,z,w,v\rightarrow1 \atop u\rightarrow y}
(\hat{\delta} s-q  y)
  D_{\frac{zu-yw}{\delta^2}}\left(\delta t+\frac{ w- y}{\delta}\right)\\[6pt]
  &\qquad=-\frac{\pi e^{\frac{y-x}{4}t^2}{_1F_1}\left(\frac{1-y}{2(x-y)};\frac{1}{2};\frac{x-y}{2}t^2\right)}{\sqrt{2}\Gamma\left(\frac{1}{2}+\frac{1-x}{2(y-x)}\right)\Gamma\left(\frac{1}{2}+\frac{1-y}{2(x-y)}\right)}
  -\frac{\pi e^{\frac{y-x}{4}t^2}t\sqrt{-2(x-y)^2}{_1F_1}\left(\frac{1}{2}+\frac{1-y}{2(x-y)};\frac{3}{2};\frac{x-y}{2}t^2\right)}{\Gamma\left(\frac{1-x}{2(y-x)}\right)\Gamma\left(\frac{1-y}{2(x-y)}\right)}\\[6pt]
  &\qquad\qquad+\frac{\pi e^{\frac{y-x}{4}t^2}\sqrt{y-x}{_1F_1}\left(\frac{1-y}{2(x-y)};\frac{1}{2};\frac{x-y}{2}t^2\right)}{\Gamma\left(\frac{1-x}{2(y-x)}\right)\Gamma\left(\frac{1}{2}+\frac{1-y}{2(x-y)}\right)}
  +\frac{\pi e^{\frac{y-x}{4}t^2}t\sqrt{x-y}{_1F_1}\left(\frac{1}{2}+\frac{1-y}{2(x-y)};\frac{3}{2};\frac{x-y}{2}t^2\right)}{\Gamma\left(\frac{1-y}{2(x-y)}\right)\Gamma\left(\frac{1}{2}+\frac{1-x}{2(y-x)}\right)}\\[6pt]
  &\qquad=-A_1-A_2+A_3+A_4,
\end{align*}
and
\begin{align*}
&\lim_{ y,z,w,v\rightarrow1 \atop u\rightarrow y}
(pw-\delta r)
 D_{\frac{xv-yw}{\hat{\delta}^2}}\left(\hat{\delta} t+\frac{ y- w}{\hat{\delta}}\right)\\[6pt]
  &\qquad=\frac{\pi e^{\frac{x-y}{4}t^2}{_1F_1}\left(\frac{1-x}{2(y-x)};\frac{1}{2};\frac{y-x}{2}t^2\right)}{\sqrt{2}\Gamma\left(\frac{1}{2}+\frac{1-x}{2(y-x)}\right)\Gamma\left(\frac{1}{2}+\frac{1-y}{2(x-y)}\right)}
  +\frac{\pi e^{\frac{x-y}{4}t^2}t\sqrt{-2(x-y)^2}{_1F_1}\left(\frac{1}{2}+\frac{1-x}{2(y-x)};\frac{3}{2};\frac{y-x}{2}t^2\right)}{\Gamma\left(\frac{1-x}{2(y-x)}\right)\Gamma\left(\frac{1-y}{2(x-y)}\right)}\\[6pt]
  &\qquad\qquad-\frac{\pi e^{\frac{x-y}{4}t^2}\sqrt{x-y}{_1F_1}\left(\frac{1-x}{2(y-x)};\frac{1}{2};\frac{y-x}{2}t^2\right)}{\Gamma\left(\frac{1-y}{2(x-y)}\right)\Gamma\left(\frac{1}{2}+\frac{1-x}{2(y-x)}\right)}
  -\frac{\pi e^{\frac{x-y}{4}t^2}t\sqrt{y-x}{_1F_1}\left(\frac{1}{2}+\frac{1-x}{2(y-x)};\frac{3}{2};\frac{y-x}{2}t^2\right)}{\Gamma\left(\frac{1-x}{2(y-x)}\right)\Gamma\left(\frac{1}{2}+\frac{1-y}{2(x-y)}\right)}\\[6pt]
  &\qquad=B_1+B_2-B_3-B_4.
\end{align*}
Since
\begin{equation}\label{FtoeF}
{_1F_1}(a;b;z)=e^z{_1F_1}(b-a;b;-z),
\end{equation}
we deduce that
\begin{align}
  _1F_1\left(\frac{1-x}{2(y-x)};\frac{1}{2};\frac{y-x}{2}t^2\right)&=e^{\frac{y-x}{2}t^2}{_1F_1}\left(\frac{1-y}{2(x-y)};\frac{1}{2};\frac{x-y}{2}t^2\right),\label{ff1}\\[6pt]
  {_1F_1}\left(\frac{1}{2}+\frac{1-x}{2(y-x)};\frac{3}{2};\frac{y-x}{2}t^2\right)&=e^{\frac{y-x}{2}t^2}{_1F_1}\left(\frac{1}{2}+\frac{1-y}{2(x-y)};\frac{3}{2};\frac{x-y}{2}t^2\right).\label{ff2}
\end{align}
Plugging \eqref{ff1} into $B_1$ and $B_3$,
and plugging \eqref{ff2} into $B_2$ and $B_4$, we see that
$A_1=B_1$, $A_2=B_2$ and
\begin{align*}
  &A_3-B_3= \pi e^{\frac{y-x}{4}t^2}\Theta{_1F_1}\left(\frac{1-y}{2(x-y)};\frac{1}{2};\frac{x-y}{2}t^2\right),\\[6pt]
  &A_4-B_4= -\pi e^{\frac{y-x}{4}t^2}\Theta t{_1F_1}\left(\frac{1}{2}+\frac{1-y}{2(x-y)};\frac{3}{2};\frac{x-y}{2}t^2\right),
\end{align*}
which implies
\begin{align}\label{denominator}
  &\lim_{ y,z,w,v\rightarrow1 \atop u\rightarrow y}
(\hat{\delta} s-q  y)
  D_{\frac{zu-yw}{\delta^2}}\left(\delta t+\frac{ w- y}{\delta}\right)+(pw-\delta r)
 D_{\frac{xv-yw}{\hat{\delta}^2}}\left(\hat{\delta} t+\frac{ y- w}{\hat{\delta}}\right)\\[6pt]
  &\qquad\quad=\pi e^{\frac{y-x}{4}t^2}\Theta \left({_1F_1}\left(\frac{1-y}{2(x-y)};\frac{1}{2};\frac{x-y}{2}t^2\right)-t{_1F_1}\left(\frac{1}{2}+\frac{1-y}{2(x-y)};\frac{3}{2};\frac{x-y}{2}t^2\right)\right).\notag
\end{align}
Therefore, combining \eqref{numerator} and \eqref{denominator} gives \eqref{genofTnxy},
which completes the proof.
\qed
%\begin{thm}\label{corofTn}
% We have
%\begin{equation}\label{genofTnxy}
%\sum_{n=0}^\infty T_n (x,y) \frac{t^n}{n!}=\frac{e^{\frac{x-y}{2} t^2 }}{{_1F_1}\left(\frac{1-y}{2 (x-y)};\frac{1}{2};\frac{x-y}{2} t^2\right)- t\,{_1F_1}\left(\frac{1}{2}+\frac{1-y}{2(x-y)};\frac{3}{2};\frac{x-y}{2} t^2\right) }.
%\end{equation}
%\end{thm}

Letting $y\rightarrow x$ in Theorem \ref{corofTn} also leads to the result of Gessel \cite[A008971]{OEIS}.

\noindent{\emph{Proof of Theorem \ref{Gesselthm}}.} One can check that
$$
\lim_{y\rightarrow x}{_1F_1}\left(\frac{1-y}{2 (x-y)};\frac{1}{2};\frac{x-y}{2} t^2\right)=\sum_{n=0}^\infty \frac{(\sqrt{1-x}t)^{2n}}{(2n)!}=\cosh(\sqrt{1-x}t)
$$
and
$$
\lim_{y \rightarrow x}{_1F_1}\left(\frac{1}{2}+\frac{1-y}{2(x-y)};\frac{3}{2};\frac{x-y}{2} t^2\right)=\frac{1}{\sqrt{1-x}t}\sum_{n=0}^\infty  \frac{(\sqrt{1-x}t)^{2n+1}}{(2n+1)!}=\frac{1}{\sqrt{1-x}t}\sinh(\sqrt{1-x}t).
$$

Thus we have
\begin{align*}
\lim_{y \rightarrow x}\sum_{n=0}^\infty T_n (x,y) \frac{t^n}{n!}&=\lim_{y \rightarrow x}\frac{e^{\frac{x-y}{2} t^2 }}{{_1F_1}\left(\frac{1-y}{2 (x-y)};\frac{1}{2};\frac{x-y}{2} t^2\right)- t\,{_1F_1}\left(\frac{1}{2}+\frac{1-y}{2(x-y)};\frac{3}{2};\frac{x-y}{2} t^2\right) }\\[6pt]
 &=\frac{\sqrt{1-x}}{\sqrt{1-x}\cosh(\sqrt{1-x}t)-\sinh(\sqrt{1-x}t)}\\[6pt]
 &=\sum_{n=0}^\infty T_n(x)\frac{t^n}{n!}.
\end{align*}
We complete the proof. \qed

Furthermore, setting $y\rightarrow 1$ and $x\rightarrow 1$ in \eqref{Tnxy}, respectively, we  deduce that
\begin{equation*}
\bar{T}_n(x):=T_n(x,1)=\sum_{k\geq0}T_1(n,k)x^k\quad\mbox{and}\quad \tilde{T}_n(y):=T_n(1,y)=\sum_{k\geq0}T_2(n,k)y^k,
\end{equation*}
where $T_1(n,k)$  is
the number of permutations on $[n]$ with $k$ exterior peaks of pattern 132,
and $T_2(n,k)$ is the number of permutations on $[n]$ with $k$ exterior peaks of pattern 231.

\begin{cor}
  We have
  \begin{equation}\label{gentx}
    \sum_{n=0}^{\infty}{\bar{T}_n(x)}\frac{t^n}{n!}=\frac{e^{\frac{x-1}{2}t^2}}{1-\int_0^te^{\frac{x-1}{2}s^2}ds}
  \end{equation}
  and
  \begin{equation}\label{genty}
    \sum_{n=0}^{\infty}\tilde{T}_n(y)\frac{t^n}{n!}=\frac{1}{1-\int_0^te^{\frac{y-1}{2}s^2}ds}.
  \end{equation}
\end{cor}
\pf
Since $$_1F_1\left(\frac{1}{2};\frac{3}{2};-z^2\right)=\frac{\sqrt{\pi}}{2z}\erf(z),$$  we observe that
\begin{equation}\label{sub1}
  \lim_{y\rightarrow 1}{_1F_1}\left(\frac{1}{2}+\frac{1-y}{2(x-y)};\frac{3}{2};\frac{x-y}{2} t^2\right)
  =\frac{\sqrt{\pi}}{\sqrt{2(1-x)}t}\erf\left(\sqrt{\frac{1-x}{2}}t\right).
\end{equation}
Then by Theorem  \ref{corofTn}, setting $y\rightarrow1$ in \eqref{genofTnxy}  and utilizing \eqref{sub1},
we obtain
\begin{equation}\label{sub11}
  \sum_{n=0}^{\infty}{\bar{T}_n(x)}\frac{t^n}{n!}=\frac{\sqrt{2(1-x)}e^{\frac{x-1}{2}t^2}}{\sqrt{2(1-x)}-\sqrt{\pi}\erf\left(\sqrt{\frac{1-x}{2}}t\right)}.
\end{equation}
With the help of \eqref{1f1toerf},
we deduce that
\begin{equation}\label{sub2}
  \lim_{x\rightarrow 1}{_1F_1}\left(\frac{1}{2}+\frac{1-y}{2(x-y)};\frac{3}{2};\frac{x-y}{2} t^2\right)
  =\frac{\sqrt{\pi}}{\sqrt{2(1-y)}t}e^{\frac{1-y}{2}t^2}\erf\left(\sqrt{\frac{1-y}{2}}t\right).
\end{equation}
Thus, setting $x\rightarrow1$ in \eqref{genofTnxy} and using \eqref{sub2}, we see that
\begin{equation}\label{sub22} \sum_{n=0}^{\infty}\tilde{T}_n(y)\frac{t^n}{n!}=\frac{\sqrt{2(1-y)}}{\sqrt{2(1-y)}-\sqrt{\pi}\erf\left(\sqrt{\frac{1-y}{2}}t\right)}.
\end{equation}
After replacement error functions by integration in \eqref{sub11} and \eqref{sub22},
we obtain \eqref{gentx} and \eqref{genty}. This completes the proof.
\qed

Setting $x=0$ in \eqref{gentx},
we obtain the generating function of the number of permutations on $[n]$ with no exterior peaks of pattern 132
\begin{equation*}
    \sum_{n=0}^{\infty}T_1(n,0)\frac{t^n}{n!}=\frac{e^{-\frac{t^2}{2}}}{1-\int_0^te^{-\frac{s^2}{2}}ds},
\end{equation*}
which recasts the formula obtained by Kiteav in \cite[Theorem 6]{Kitaev-2003}.
Since $\pi_0\pi_1\pi_2$ cannot be a exterior peak of pattern 231,
$T_2(n,k)$ is also the number of  permutations on $[n]$ with $k$ consecutive pattern $231$.
Thus, \eqref{genty} reproduces the result of Elizalde-Noy given in \cite[Theorem 4.1]{Elizalde-Noy 2003}.
Similarly, setting $y=0$ in \eqref{genty} yields the generating functions of permutations on $[n]$
avoiding the consecutive patterns $231$:
\begin{equation*}
    \sum_{n=0}^{\infty}T_2(n,0)\frac{t^n}{n!}=\frac{1}{1-\int_0^te^{-\frac{s^2}{2}}ds},
\end{equation*}
which is also a special case of \cite[Theorem 4.1]{Elizalde-Noy 2003} and \cite[Theorem 12]{Kitaev-2005}.

As the conclusion of this section,
by combining Theorem \ref{thmofGen z} and Theorem \ref{thmofGen v},
we give the joint distribution of exterior peaks
of pattern 132 and of pattern 231 over
alternating permutations on $[n]$.
As defined in \cite{Stanley-2012},
we say that a permutation $\pi=\pi_1\pi_2\cdots\pi_n\in\mathfrak{S}_n$ is alternating (or zigzag or down-up)
if $\pi_1>\pi_2<\pi_3>\pi_4<\cdots$.
Recall that $T^A(n,i,j)$ is the number of alternating  permutations on $[n]$ with
$i$ exterior peaks of pattern 132 and $j$ exterior peaks of pattern 231, and the generating function
$$
T^A_n(x,y)=\sum_{i,j \geq 0} T^A(n,i,j)x^i y^j.
$$

From \eqref{TA}, we see that
\begin{equation*}
  \lim_{x,y\rightarrow1}{_1F_1}\left(\frac{x}{2(x-y)};\frac{3}{2};-\frac{x-y}{2}t^2\right)=\sum_{n=0}^{\infty}\frac{(-1)^n}{(2n+1)!}t^{2n}=\frac{\sin t}{t}
\end{equation*}
and
\begin{equation*}
  \lim_{x,y\rightarrow1}{_1}F_{1}\left(-\frac{y}{2(x-y)};\frac{1}{2};\frac{x-y}{2}t^2\right)=\sum_{n=0}^{\infty}\frac{(-1)^n}{(2n)!}t^{2n}=\cos t.
\end{equation*}
Thus if we set $x,y\rightarrow 1$ in \eqref{TA}, then \eqref{TA} reduces to the generating function
\begin{equation*}
  \sum_{n=0}^{\infty}E_n\frac{t^n}{n!}=\sec t+\tan t,
\end{equation*}
where the
Euler number $E_n$ is the number of alternating permutations on $[n]$.

\noindent{\emph{Proof of Theorem \ref{coroofTA}.}}
The proof is divided  into the following cases depending on the parity of $n$.
For the case that $n$ is even, we let $z,v\rightarrow1$ and $y,w\rightarrow0$ and $u\rightarrow y$ in \eqref{pnijk}.
One can check that the permutations not vanished are alternating permutations with even length. Thus we have
\begin{equation*}
    P_{2k}(x,0,1,0,y,1)=T^A_{2k}(x,y)  \quad\mbox{and}\quad P_{2k+1}(x,0,1,0,y,1)=0,
\end{equation*}
where $k$ is a nonnegative integer.
By setting $z,v\rightarrow1$ and $y,w\rightarrow0$ and $u\rightarrow y$ in Theorem \ref{thmofGen z},
using formula \eqref{FtoeF}, we obtain that
\begin{equation}\label{TAeven}
\sum_{k=0}^\infty T^A_{2k} (x,y) \frac{t^{2k}}{(2k)!}=\frac{e^{\frac{x-y}{2}t^2}}{{_1}F_{1}\left(-\frac{y}{2(x-y)};\frac{1}{2};\frac{x-y}{2}t^2\right)}.
\end{equation}

For the case that {$n\geq3$} is odd, we set $z,v\rightarrow1$ and $y,w\rightarrow0$ and $u\rightarrow y$ in \eqref{qnijkl}.
Notice that only the alternating permutations with odd length are not vanished in this case. By the labeling rules in proof of Theorem \ref{thmofQ}, we have
\begin{equation*}
   Q_{2k+1}(x,0,1,0,y,1)/y=T^A_{2k+1}(x,y)
   \quad\mbox{and}\quad Q_{2k}(x,0,1,0,y,1)/y=0,
\end{equation*}
where $k$ is a positive integer.
By letting $z,v\rightarrow1$ and $y,w\rightarrow0$ and $u\rightarrow y$ in Theorem \ref{thmofGen v}, then using formula
$$(1+a-b){_1F_1}(a;b;z)=a{_1F_1}(a+1;b;z)+(1-b){_1F_1}(a;b-1;z),$$  we derive that
\begin{equation}\label{TAodd}
  \sum_{k=1}^{\infty} T^A_{2k+1}(x,y)\frac{t^{2k+1}}{(2k+1)!}=
  \frac{te^{\frac{x-y}{2}t^2}{_1F_1}\left(\frac{x}{2(x-y)};\frac{3}{2};-\frac{x-y}{2}t^2\right)}
  {{_1F_1}\left(-\frac{y}{2(x-y)};\frac{1}{2};\frac{x-y}{2} t^2\right)}-1.
\end{equation}
Therefore, summing $T^A_1(x,y)=1$ and \eqref{TAeven}, \eqref{TAodd} together yields \eqref{TA},
which completes the proof.
\qed

\vspace{0.5cm}

\noindent{\bf Acknowledgements}.
The authors appreciate Professor William Y.C. Chen for the valuable suggestions.

\end{document}